%% file: paper.tex
\documentclass{article}
\usepackage{amsfonts,amsmath,amssymb}
\usepackage{algorithm,algpseudocode}
\usepackage{color}
\usepackage{graphicx}
\usepackage{url}
\usepackage{verbatim}

\DeclareMathOperator{\argmin}{argmin}

\DeclareMathOperator{\Span}{span}

\newcommand{\Range}[1]{\ensuremath{\mathcal{R}\left( #1 \right)}}

\newcommand{\thinColon}{%
  \nobreak
  \mskip1mu 
  \mathpunct{}%
  \nonscript
  \mkern-\thinmuskip
  {:}%
  \mskip1mu
  \relax
}
\newcommand{\IndexRange}[2]{\ensuremath{{#1}\thinColon{}{#2}}}

\begin{document}

\title{Fault-tolerant linear solvers via selective reliablity}
\author{Patrick G.\ Bridges, Kurt B.\ Ferreira, Michael A.\ Heroux, and Mark Hoemmen}
% Department of Computer Science
% University of New Mexico
% Albuquerque, NM 87131
% Email: bridges@cs.unm.edu
%
% Sandia National Laboratories
% P.O.\ Box 5800, MS 1320 
% Albuquerque, NM 87185-1320, USA
% Email: \{kbferre,maherou,mhoemme\}@sandia.gov}
\maketitle

\input{abstract}
\input{intro}

\input{related-work}
\input{fault-char}
\input{models}
\input{desired-properties}
\input{alg}
\input{programming-model}
\input{numerical-experiments}

\input{interface}
\input{impl}

\input{experiments}
\input{conclusion}
\input{ack}

\bibliographystyle{siam}
\bibliography{paper}
\end{document}

%% file: abstract.tex
Energy increasingly constrains modern computer hardware, yet
protecting computations and data against errors costs energy.  This
holds at all scales, but especially for the largest parallel computers
being built and planned today.  As processor counts continue to grow,
the cost of ensuring reliability consistently throughout an
application will become unbearable.  However, many algorithms only
need reliability for certain data and phases of computation.  This
suggests an algorithm and system codesign approach.  We show that if
the system lets applications apply reliability selectively, we can
develop algorithms that compute the right answer despite faults.
These ``fault-tolerant'' iterative methods either converge eventually,
at a rate that degrades gracefully with increased fault rate, or
return a clear failure indication in the rare case that they cannot
converge.  Furthermore, they store most of their data unreliably, and
spend most of their time in unreliable mode.

We demonstrate this for the specific case of detected but
uncorrectable memory faults, which we argue are representative of all
kinds of faults.  We developed a cross-layer application / operating
system framework that intercepts and reports uncorrectable memory
faults to the application, rather than killing the application, as
current operating systems do.  The application in turn can mark memory
allocations as subject to such faults.  Using this framework, we wrote
a fault-tolerant iterative linear solver using components from the
Trilinos solvers library.  Our solver exploits hybrid parallelism (MPI
and threads).  It performs just as well as other solvers if no faults
occur, and converges where other solvers do not in the presence of
faults.  We show convergence results for representative test problems.
Near-term future work will include performance tests.

%% file: intro.tex
%%%%%%%%%%%%%%%%%%%%%%%%%%%%%%%%%%%%%%%%%%%%%%%%%%%%%%%%%%%%%%%%%%%%%%
\section{Introduction}\label{S:intro}
%%%%%%%%%%%%%%%%%%%%%%%%%%%%%%%%%%%%%%%%%%%%%%%%%%%%%%%%%%%%%%%%%%%%%%

Computational scientists tend to think of computer systems as reliable
digital devices.  Decades of experience confirmed this view, because
any faults that did occur were infrequent enough that hardware or
system software fault detection and correction schemes could handle
them.  However, many system designers predict that reliability will
decline on future computers, especially for very high-end computers
built of millions of
components~\cite{Miskov-Zivanov:2007:SER:1266366.1266680,Karnik:2004:CSE:1032295.1032595}.
This is because current and future hardware is energy constrained.
All hardware or software methods for improving reliability require
energy, because they all involve redundant storage and computation.
This includes redundant data encoding (such as Reed-Solomon codes),
and redundant arithmetic computation in space or time.  Extreme-scale
hardware is particularly energy constrained, and its large number of
components makes the failure of any one of them more likely,
increasing the demands on fault detection and correction.  There are
many efforts in the hardware development community to understand these
issues, for example~\cite{Chishti:2009:ICL:1669112.1669126,
  Yang:2005:LRS:1120725.1120957}.  Some studies already indicate that
faults are appearing at the user
level~\cite{Haque:2010:HDS:1844765.1845231}. However, without fault
detection in the user code, these faults are not always noticed, even
though they may lead to incorrect results.

Most existing approaches to fault-tolerant algorithm development
assume that a fault can occur at any time during program execution.
In this paper we explore the use of variable reliability to develop
algorithms that perform most computations using a less reliable
computing mode, but perform some computations in a special, more
highly-reliable environment.  Using this approach, we show that with
modest modifications, common iterative methods can exhibit reliable
behavior even if faults occur during the computation.  Furthermore, we
believe this basic approach can be applied to many classes of
algorithms such that, by performing a small fraction of an algorithm's
computations in highly-reliable mode, we can continue to make progress
in our computations in the presence of some system unreliability.

Both hardware and system software architects must take ever more
extreme measures to maintain the illusion of reliability with
increasingly unreliable hardware.  Yet, many algorithms do not need
this level of reliability everywhere.  Reducing energy requirements
for future computers requires an algorithm / system codesign approach.
We are using this research as a model to improve collaboration between
these two fields.

%% file: related-work.tex
%%%%%%%%%%%%%%%%%%%%%%%%%%%%%%%%%%%%%%%%%%%%%%%%%%%%%%%%%%%%%%%%%%%%%%
\section{Related work}\label{S:related}
%%%%%%%%%%%%%%%%%%%%%%%%%%%%%%%%%%%%%%%%%%%%%%%%%%%%%%%%%%%%%%%%%%%%%%

% TODO (mfh 24 Mar 2011) Make clear that we're talking about
% fault-tolerant _numerical_ algorithms.  Fault-tolerant algorithms in
% general tend to relate to things like peer-to-peer storage and
% computing, and distributed decision problems.

Fault-tolerant algorithms have long been a topic of research.  In
numerical linear algebra, most fall within the category of
\emph{algorithm-based fault tolerance} (ABFT) (see e.g.,
\cite{huang1984algorithm}).  Such approaches are interesting research,
but often do not fully address the needs of applications.  In
particular, ABFT methods attempt to detect faults during the execution
of some function such as a solver, and then recover solver state via
metadata collected during execution or basic mathematical properties
known about the algorithm.  However, such approaches are impractical
since solver state is only one portion of the total application state.
If application state is not also recovered, the solver state is
irrelevant.  Furthermore, solver state is easily regenerated if
application state is recovered.  As a result, ABFT methods are not
presently used in applications as far as we know.  ABFT methods can
become relevant if we can finally have in place the vertically
integrated resilience capabilities mentioned in the context of hard
fault situations.  In this situation, faults detected and resolved in
the solver can remain relevant if the application has also managed to
recover its corresponding state.

Other authors have empirically investigated the behavior of iterative
solvers when soft faults occur (e.g.,
\cite{Bronevetsky:2008:Soft,howle2010soft:Copper}), developed an
approximate restart scheme for recovering from the loss of a node's
data \cite{Langou:2007:RPI:1350656.1350657}, or even developed more
energy-conserving hardware cache error correction schemes, based on
observations of iterative methods' cache use
\cite{malkowski2010analyzing}.  ``Asynchronous'' or ``chaotic''
iterations (see e.g., \cite{bahi2007parallel} for a bibliography) are
linear solvers designed to tolerate message delays when applying the
matrix in parallel, for certain classes of matrices.  However, as far
as we know, no one has yet developed iterative solver algorithms
specifically to handle soft faults in computations and data.

%%% Local Variables: 
%%% mode: latex
%%% TeX-master: "paper"
%%% End: 

%% file: fault-char.tex
%%%%%%%%%%%%%%%%%%%%%%%%%%%%%%%%%%%%%%%%%%%%%%%%%%%%%%%%%%%%%%%%%%%%%%
\section{Fault characterization}\label{S:faultchar}
%%%%%%%%%%%%%%%%%%%%%%%%%%%%%%%%%%%%%%%%%%%%%%%%%%%%%%%%%%%%%%%%%%%%%%

In this paper, we use \emph{fault} to mean an abnormal operating
condition of the computer system, which affects a running routine (in
this case a linear solver) in some way.  The routine \emph{fails} only
when one or more faults causes it to compute the wrong answer.  That
is, faults occur inside a routine; failure refers to the routine's
output, which does not meet the caller's success criteria.  This
distinction between faults and failures is a simplified version of the
multilevel model of software reliability presented in
\cite{parhami1997defect}.  This definition nests: for example, if a
nonlinear solver calls a linear solver repeatedly, the linear solver
may produce a solution with residual norm greater than the caller's
tolerance (i.e., ``fail'') on occasion, but the nonlinear solver may
still converge.  Thus, failure from the linear solver's perspective
may be a fault but not a failure from the nonlinear solver's
perspective.  The rest of this paper considers faults and failures
from the linear solver's perspective.  We leave studies of algorithms
that consume linear solvers' output (such as nonlinear solvers,
optimization algorithms, and implicit methods for solving
time-dependent systems of ordinary differential equations) for future
work.  
% Howle \cite{howle2010soft:Copper} discussed empirical observations
% of nonlinear solver behavior under linear solver failures.

In this paper, we give two classifications of faults.  The first is
\emph{hard} vs.\ \emph{soft}:
\begin{itemize}
\item \textbf{Hard faults:} Cause program interruption and are outside
  the scope of what the executable program can directly detect.  These
  faults can result from hardware failure or from data integrity
  faults that lead to an incorrect execution path.  
\item \textbf{Soft faults:} Do not cause immediate program
  interruption and are detectable via introspection by user code.
  Soft faults occur as ``bit flips'' such as incorrect floating point
  or integer data, or perhaps incorrect address values that still
  point to valid user data space.  Although it is difficult to detect
  all soft faults, some modest amount of introspection can be very
  effective at dramatically reducing their impact.
\end{itemize}
An example of a hard fault would be the operating system crashing,
causing the program to stop executing.  (This would not be a failure
if the system then restarts the program from a checkpoint, and the
program completes and produces the correct answer.)  In our
experience, detecting and recovering from hard faults requires a
concerted effort from all levels of the hardware and software stack.
Although there may be algorithmic research required for this effort,
the primary need is to determine roles, responsibilities and protocols
for communicating between layers.  This activity is underway in some
layers, but is only starting to be addressed in a comprehensive way.

\begin{figure}[!h]
\begin{center}
\includegraphics[scale=0.85]{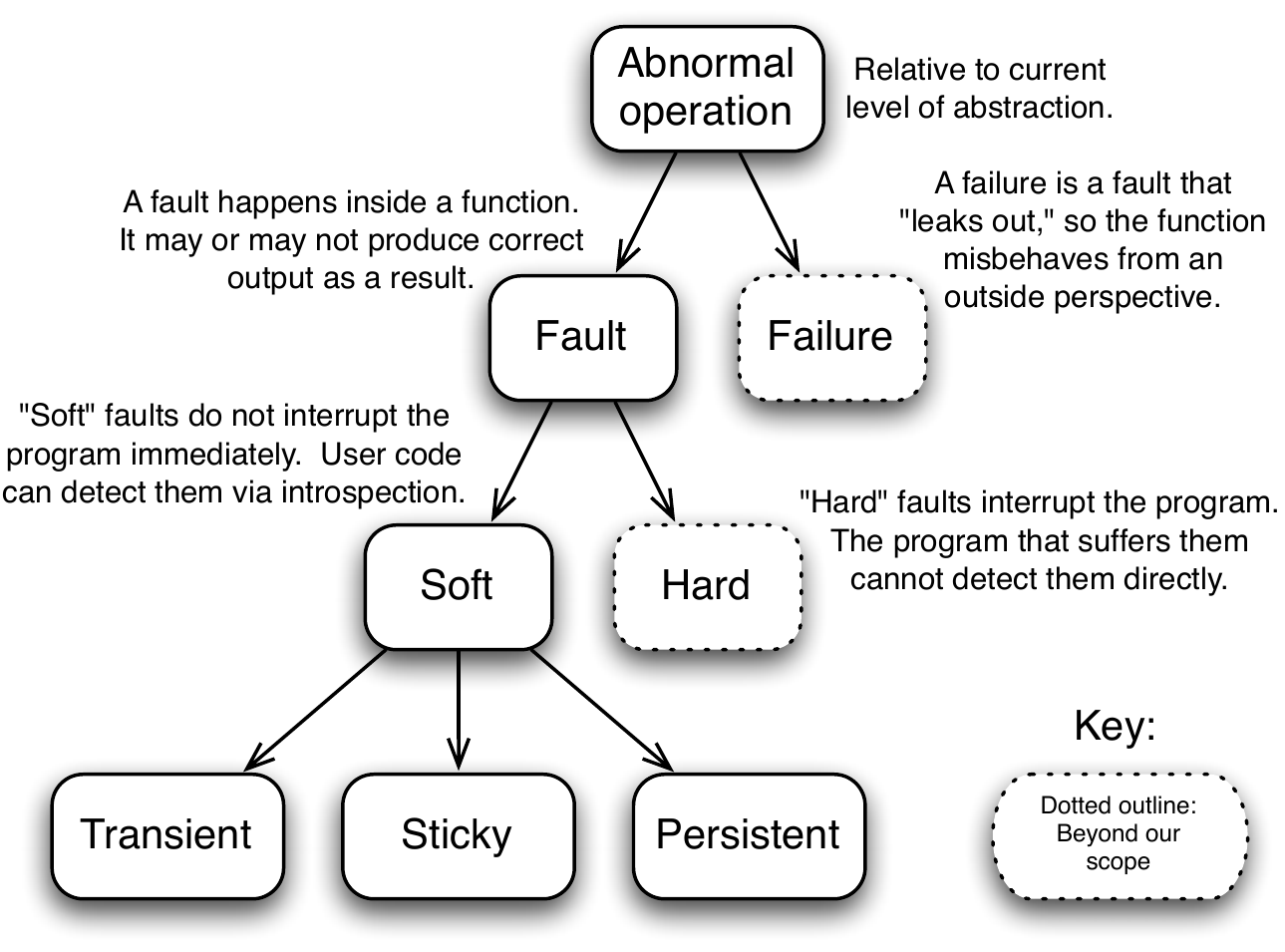}
\end{center}
\caption{\label{fig:FaultClassification}Classification of faults.
  Hard faults are outside the scope of our effort.  We address soft
  faults in several ways.}
\end{figure}

%
% mfh 24 Mar 2011: "Persistent" or "permanent"?  I guess "persistent,"
% since nothing is really permanent (buy a new computer and resolve
% the problem from the beginning, in the worst case).
%
The second characterization applies only to soft faults, and describes
their temporal behavior:
\begin{itemize}
\item \textbf{Persistent fault:} The incorrect bit pattern will not
  change as execution proceeds.  Example: The primary source of a data
  value (and any subsequent copies) are incorrect, so there is no
  ability to restore correct state.
\item \textbf{Sticky fault:} The incorrect bit pattern can be
  corrected by direct action.  Example: A backup source for the data
  exists and can be used to restore correct state.
\item \textbf{Transient fault:} The incorrect pattern occurs
  temporarily. Example: Data in a cache is incorrect, but the correct
  value is still present in main memory and the cache value is
  flushed.
\end{itemize}
Figure \ref{fig:FaultClassification} illustrates the relationship
between the two characterizations of faults.

%%%%%%%%%%%%%%%%%%%%%%%%%%%%%%%%%%%%%%%%%%%%%%%%%%%%%%%%%%%%%%%%%%%%%%
\subsection{Potential for Soft Fault Detection and Correction}
\label{subsect:Potential}
\label{SS:faultchar:potential} % Support both notations for now
%%%%%%%%%%%%%%%%%%%%%%%%%%%%%%%%%%%%%%%%%%%%%%%%%%%%%%%%%%%%%%%%%%%%%%

Although recovering from hard faults requires a coordinated effort
across software and hardware layers, at least some soft faults can be
effectively detected and corrected by user code.  Furthermore,
practically speaking, some applications spend much of their
computation time in a small portion of the total program lines of
code.  Such applications can benefit from introducing fault-oriented
introspection into that portion of the software.  This situation
occurs frequently in applications that generate and solve large linear
systems of equations.  In many cases, 80\% or more of the computation
time is spent in the linear solver. As problem sizes and processor
counts increase, the solver can take more than 99\% of the total
execution time~\cite{Lin:2010:TLM:1838773.1839098}.  If we can
incorporate introspection into the solvers for these cases, we can
dramatically reduce the impact of soft faults.

One of the challenges for future system designers is determining how
much fault resilience should be designed into the system.
Historically, hardware and system architects have been very aggressive
in capturing faults, so much so that users rarely experience a fault
during normal system use.  In the future, such approaches may be too
expensive, resulting in a default reliability that must always be
scrutinized.  With this in mind, we introduce the concept of high vs.\
bulk reliability:
\begin{itemize}
\item \textbf{Bulk reliability:} The default reliability exhibited by
  system in normal execution mode.  As system feature sizes shrink and
  component counts increase, we expect that bulk reliability will
  decrease to the point where users will need to pay attention to
  potential errors.
\item \textbf{High reliability:} A special, presumably
  software-enabled mode, such that the user can declare data storage
  regions, data paths and execution regions that have better than bulk
  reliability.
\end{itemize}

Presently most algorithms lack robustness in the presence of soft
faults.  A single soft fault will not be detected and will eventually
result in catastrophic failure.  Assuming we have high reliability
mechanisms in future programming environments, we have new
opportunities for redesigning algorithms.  Specifically, we seek
algorithm designs such that decay in progress is proportional to the
number of soft faults, at least in practice.

In this paper we focus on preconditioned iterative methods, and
particularly on variants of GMRES (the Generalized Minimal Residual
method \cite{saad1986gmres}).  We do so because, as we mentioned, many
applications spend the vast majority of execution time in the solver,
and GMRES is one of the most robust and popular methods for
challenging problems.  However, the approach we use is applicable to
many algorithms.  In fact, we believe that most, and maybe all,
algorithms can eventually have fault-resilient formulations that
introduce a very small runtime overhead while practically achieving
the convergence equivalent to doing all computations in high
reliability mode.

%%% Local Variables: 
%%% mode: latex
%%% TeX-master: "paper"
%%% End: 

%% file: models.tex
%%%%%%%%%%%%%%%%%%%%%%%%%%%%%%%%%%%%%%%%%%%%%%%%%%%%%%%%%%%%%%%%%%%%%%
\section{Models of reliability}\label{S:reliability-models}
%%%%%%%%%%%%%%%%%%%%%%%%%%%%%%%%%%%%%%%%%%%%%%%%%%%%%%%%%%%%%%%%%%%%%%

In this section, we describe models of reliability that fault-tolerant
numerical algorithms could use.  The main goal of these models is to
help algorithm developers reason about the quality of the computed
solution.  Without the promise of reliability for selected data and
computations, no implementation of an algorithm can promise anything
about the final result.  Thus, all the models we propose in this
section allow programmers to demand reliability as needed, and to
allow data and control to flow between reliable and unreliable parts
of the program.

A second goal of our reliability models is to convert hard faults into
soft faults whenever our algorithms can handle the latter effectively.
Reliability models govern the distinction between hard and soft
faults.  For example, the \emph{fail-stop} model ensures that either
the data and computations are reliable, or the program terminates with
minimal side effects; it tries to turn all soft faults into hard
faults.  Current numerical algorithms assume a fail-stop model, which
we assert can be relaxed in many cases.  As long as algorithms can
deal with soft faults without a large time-to-solution penalty,
reducing the number of hard faults will improve performance by
avoiding restarts and allowing reduction of the checkpoint frequency.
It may even improve reliability, for example by avoiding the
catastrophic situation of a second hard fault during recovery from one
hard fault.

We begin in Section \ref{SS:reliability-models:UQ} by asking whether
statistics could help us avoid considering models of reliability, and
showing that it does not.  Section \ref{SS:reliability-models:sandbox}
describes the ``sandbox'' model, which is the most general reliability
model our fault-tolerant algorithms can use.  The algorithm presented
in Section \ref{S:alg} can work even in this model, but finer-grained
models will allow us to define its convergence behavior more
precisely.  Therefore, we conclude with some desired features of
finer-grained models in Section \ref{SS:reliability-models:features}.

% , and conclude in Section \ref{SS:reliability-models:finer} with
% some proposals for finer-grained models.

%%%%%%%%%%%%%%%%%%%%%%%%%%%%%%%%%%%%%%%%%%%%%%%%%%%%%%%%%%%%%%%%%%%%%%
\subsection{Statistical ``model''}\label{SS:reliability-models:UQ}
%%%%%%%%%%%%%%%%%%%%%%%%%%%%%%%%%%%%%%%%%%%%%%%%%%%%%%%%%%%%%%%%%%%%%%

Increasingly numerical simulations use statistical techniques to
account for uncertainty in the data as well as in the mathematical
model.  Many people refer to the study of representing and quantifying
such uncertainties as \emph{uncertainty quantification} (UQ).  It
seems reasonable that we also could apply these techniques to account
for possibly unreliable solves, that is, ``roll up'' the solver's
uncertainty in that of the application itself.  This would not require
new solver algorithms or implementations.  Instead, the problem would
be solved multiple times using existing solvers, and statistics would
be used to remove ``outliers'' and identify the most ``believable''
solution.  This would comprise a ``model'' of reliability based on
statistical belief, rather than on any guarantees made by the system
or solver.

%
% Epistemic uncertainty -- assume that the software has no
% introspection capability; we're modeling something which is not
% inherently stochastic (``aleatory uncertainty''), but rather
% something which we aren't able to measure directly.  Actually
% no, says Brian Adams, we have enough data (from Google, etc. --
% ``DRAM Errors in the Wild'') to model it probabilistically,
% therefore it's aleatory (a probabilistic model).
%
% Could still use statistics with the inner-outer iteration to model
% belief in the quality of the inner solve...
%

This ``model'' is no model of reliability at all.  It implicitly
assumes that faults may only occur in the solver, and that the
statistical analysis that identifies the most believable solution is
free of faults.  In fact, these assumptions define the ``sandbox''
model of reliability described in the next section
(\ref{SS:reliability-models:sandbox}).  Nevertheless, one might
consider using statistical analysis to improve fault tolerance, in
combination with a satisfactory fault model.  We do not think this
should be applied na\"ively to existing fault-intolerant solvers, for
two reasons.  First, it may require running many solves to get
statistical confidence.  Second, it would throw away what numerical
analysts have learned about how iterative solvers respond to certain
kinds of faults.  For example, perturbing the matrix $A$ affects
convergence of iterative solvers more in earlier iterations than in
later iterations (see Section \ref{SS:alg:inexact}).  Finally, we will
show in this paper that iterative methods can be modified to tolerate
some soft faults, for much less cost than running a fault-intolerant
solver many times.  We do not dismiss statistical approaches
completely, though.  In particular, they may be useful to enhance
detection of faults when invoking a solver.  As we discuss in Section
\ref{SS:model:recovery}, our fault-tolerant inner-outer iteration can
save some fault recovery work if it can detect faults reliably in the
inner solves.

\begin{comment}
%
% mfh 23 Mar 2011: could compress the paragraph below
%
We assert above that computing a statistically believable solution
from a fault-intolerant solver would require prohibitively many
solves.  Solves of very large problems may take long enough, and
faults may be frequent enough on very large parallel systems, that a
fault will almost certainly occur during every solve.  If those are
mostly soft faults, then sufficient confidence about the most
believable solution may require a large number of solves.  If most are
hard faults that cause program interruption, then the necessary
checkpointing and restarting from checkpoints (the usual technique for
recovering from hard faults) would drive up the cost of each solve.
In fact, there is significant concern that without more sophisticated
fault recovery techniques, applications will spend all of their time
checkpointing, restarting from checkpoints, and even recovering from
failures experienced \emph{during} the restart.  See for example
Cappello et al.\ \cite{cappello2009exascale}.  In either case, solving
the problem enough times to identify a believable solution may be too
expensive, no matter where the line is drawn between hard and soft
faults.
\end{comment}

%%%%%%%%%%%%%%%%%%%%%%%%%%%%%%%%%%%%%%%%%%%%%%%%%%%%%%%%%%%%%%%%%%%%%%
\subsection{Sandbox model}\label{SS:reliability-models:sandbox}
%%%%%%%%%%%%%%%%%%%%%%%%%%%%%%%%%%%%%%%%%%%%%%%%%%%%%%%%%%%%%%%%%%%%%%

Relaxing reliability of \emph{all} data and computations may result in
all manner of undesirable and unpredictable behavior.  If
instructions, pointers, array indices, and boolean values used for
decisions may change arbitrarily at any time, we cannot assert
anything about the results of a computation or the side effects of the
program, even if it runs to completion without abnormal termination.
%
%One node may even corrupt the state of other nodes, for
%example by overwriting parts of memory owned by a user-space
%communication library, or performing incorrect output to a shared
%filesystem.
%
%This is why we made the distinction in Section
%\ref{SS:faultchar:potential} between ``hard'' and ``soft'' faults.
%
The least we can do is force the fault-susceptible program to execute
in a \emph{sandbox}.  This is a general idea from computer security,
that allows the execution of untrusted ``guest'' code in a partition
of the computer's state (the ``sandbox'') that protects the rest of
the computer (the ``host'') from the guest's possibly bad behavior.
Sandboxing can even protect the host against malicious code that aims
to corrupt the system's state, so it can certainly handle code subject
to unintentional faults in data and instructions.  

Sandboxes \emph{ensure isolation} of a possibly unreliable phase of
execution.  They \emph{allow data to flow between reliable and
  unreliable} phases of execution.  Also, they let the host
\emph{force guest code to stop} within a predefined finite time, or if
the host suspects the guest may have wandered astray.  This feature is
especially important in distributed-memory computation for preventing
deadlock and other failures due to ``crashed'' or unresponsive nodes.
In general, sandboxing converts some kinds of hard faults into soft
faults, and limits the scope of soft faults to the guest subprogram.

% Finite time implies that the guest obeys the bulk-synchronous model
% of parallel computation? -- No, it doesn't. It actually _forces_
% the BSP model, in some sense, but the guest can still use whatever
% model of parallel computation it likes -- asynchronous even.

Sandboxing may be implemented in different ways.
%, depending on the
%possible side effects and failure modes of the guest program, and the
%performance requirements.
For example, the guest may run in a virtual machine on the same
hardware as the host.  (See Smith and Nair
\cite{smith2005architecture} or Rosenblum
\cite{rosenblum2004reincarnation} for accessible overviews of past and
recent virtual machine technology.)  Alternately, the guest may even
run on separate hardware from the host program.  For example, guests
may run on a fast but unreliable subsystem, and the controlling host
program may run on a reliable but slower subsystem.

Here is an example of the sandbox model in operation.  In this
example, the guest program is responsible for computing sparse
matrix-vector products.  It receives a vector $x$ from the host,
computes $y := A \cdot x$ (where $A$ is the sparse matrix), and
returns $y$ to the host.  The vectors $x$ and $y$ on the host are
stored and computed with reliably.  The guest makes no promises about
the correctness of the values in the vector $y$ it returns.  It may
even return different values for the same $x$ input each time it is
invokes.  However, the sandbox ensures that the guest returns in
finite time.  (For example, it may kill the guest process if it takes
too long, and return some arbitrary solution vector if the guest did
not complete its computation.)

The fault-tolerant inner-outer iteration we will describe in Section
\ref{S:alg} uses the sandbox model.  There, the guest program performs
the task ``Solve a given linear system.''  The host program invokes
the guest repeatedly for different right-hand sides, and the host
performs its own calculations reliably.  See that section for details.
Finer-grained models of reliability may improve accuracy of the inner
solves, so we now go on to describe some desired features of these
models.

%%%%%%%%%%%%%%%%%%%%%%%%%%%%%%%%%%%%%%%%%%%%%%%%%%%%%%%%%%%%%%%%%%%%%%
\subsection{Desired features of finer-grained models}
\label{SS:reliability-models:features}
%%%%%%%%%%%%%%%%%%%%%%%%%%%%%%%%%%%%%%%%%%%%%%%%%%%%%%%%%%%%%%%%%%%%%%

The sandbox model of reliability makes only two promises of the
unreliable guest: it returns something (which may not be correct), and
it completes in fixed time.  These already suffice to construct a
working fault-tolerant iterative method, as we will show in Section
\ref{S:alg}.  However, detecting faults or being able to limit how
faults may occur would also be useful.  All of these are more
sophisticated forms of \emph{introspection}.  These finer-grained
models of reliability can be used to improve accuracy of the iterative
method, or to prove more specific promises about its convergence.  We
describe some of these below.

%%%%%%%%%%%%%%%%%%%%%%%%%%%%%%%%%%%%%%%%%%%%%%%%%%%%%%%%%%%%%%%%%%%%%%
\subsubsection{Detection}
%%%%%%%%%%%%%%%%%%%%%%%%%%%%%%%%%%%%%%%%%%%%%%%%%%%%%%%%%%%%%%%%%%%%%%

Knowing that no faults occurred in a bulk-reliability phase of
execution can avoid robustness and recovery effort in the highly
reliable phase.  We discuss this more in Section \ref{S:alg} in the
context of our inner-outer iteration.  In general, if we know that the
potentially unreliable inner solver experienced no faults, we know
that its computed intermediate state (e.g., the Krylov subspace basis)
is correct.  We can then safely use that state to accelerate the next
invocation of the inner solver.  Fault detection is therefore a
valuable feature of a reliability model, even without fault recovery.
Many error-correcting storage schemes, such as those in DRAM memory,
caches, and redundant disk storage, can detect more kinds of errors
than those which they can correct.  Extending those storage schemes to
be able to correct those additional detectable errors requires
additional hardware, energy consumption, and computation.  Thus, if
algorithms can exploit fault detection to handle faults efficiently,
they can relieve hardware of the burden of recovery.

%%%%%%%%%%%%%%%%%%%%%%%%%%%%%%%%%%%%%%%%%%%%%%%%%%%%%%%%%%%%%%%%%%%%%%
\subsubsection{Transience}
%%%%%%%%%%%%%%%%%%%%%%%%%%%%%%%%%%%%%%%%%%%%%%%%%%%%%%%%%%%%%%%%%%%%%%

Faults should look as \emph{transient} as possible.  For example,
consider solving the sparse linear system $Ax=b$ iteratively.  If
faults in the entries of $A$ persist throughout the iterative method,
the method will be solving the wrong linear system $\tilde{A}x=b$.
Worse yet, the algorithm will report that the computed approximate
solution $\tilde{x}$ has a small residual norm $\|b -
\tilde{A}\tilde{x}\|$, even though $\tilde{x}$ may be far from the
actual solution.  In contrast, many iterative methods naturally
tolerate some kinds of occasional transient faults, so unreliable
computations with only transient faults can still be useful.  Indeed,
before reliable electronic computers, the only ``computers'' were
unreliable human beings.  They could nevertheless solve real-world
problems, because human faults are usually transient.  (This is why,
when balancing a checkbook by hand, it helps to repeat the process
until one gets the same result more than once.)

% Human ``computers'' favored iterative algorithms (like Jacobi
% iteration) over Gaussian elimination for solving linear systems,
% since iterations tolerate occasional mistakes, but Gaussian
% elimination does not.

Many hardware faults are not transient.  This is particularly true of
DRAM memory faults, as described for example in Schroeder et al.\
\cite{schroeder2009dram}.  Permanent faults (which Schroeder et al.\
call ``hard errors'') due to hardware failures are much more common
than temporary faults.
%
% {\bf MAH: I wonder if we actually see these hard errors within the
%  scope of executing an application.  The impression I have is that
%  soft errors are more often ``seen'' by the app user than hard
%  errors.  We should refine this comment.}.
%
The so-called ``chip-kill'' DRAM error-correcting code (see Asanovic
et al.\ \cite{asanovic2006landscape}) was designed for the common case
of an entire DRAM module failing permanently and producing incorrect
values.  In many cases, permanent faults interrupt a running program
or even make the node fail, and are thus beyond the ability of an
application to detect.  That is, they are ``hard faults'' (see Section
\ref{S:faultchar}).  However, applications may be able to detect and
respond to these malfunctions as they first begin.  Furthermore,
``temporary'' single-bit faults may persist and accumulate into
multiple-bit faults, which some error-correcting codes cannot correct.
Eliminating correctible faults before they become uncorrectible
requires special measures (a ``memory scrubber'') that may increase
energy consumption and reduce available memory bandwidth.

This means the implementation of the reliability model likely will
have to do extra work to give the appearance of transience.  In terms
of Section \ref{S:faultchar}, the implementation must turn
``persistent'' faults into ``sticky'' or ``transient'' faults.
%
% Caches already have this feature; they refresh bad data from DRAM.
% 
For example, unreliable memory storing the sparse matrix $A$ could be
refreshed every few iterations from a reliable backing store.
Physical memory pages showing incorrect values during the refresh may
be retired and replaced with other physical pages.  The reliable
backing store approach is also useful for checkpointing, and could be
implemented with fast local storage (like flash memory).  

\begin{comment}
Some means of user annotation may be helpful in order to avoid saving
data to the reliable backing store more times than necessary.  This is
in part because iterative methods often reuse arrays in memory.  For
example, the same storage for the initial residual vector might be
reused and modified as the vector is changed for each restart cycle,
but the sparse matrix $A$ usually is not modified by the iterative
method.  Thus, the sparse matrix $A$ would only need to be saved to
reliable backing store once for all restart cycles, whereas the
initial residual vector would need to be saved once per restart cycle.
A ``transactional'' model of updates could be used to track the need
for saving data to reliable backing store: users specifically have to
acquire and release views of important arrays of data in order to
modify their contents, and the ``release'' operation indicates that
the array contents have changed and should be saved (perhaps
asynchronously) to reliable backing store.  This model is useful for
other reasons (such as execution of computational kernels on
accelerator devices like graphics processing units (GPUs)) and would
fit naturally into an already transactional programming interface such
as Kokkos (see e.g., Baker et al.\ \cite{baker2010lightweight}).
\end{comment}

% If the basis vectors
% are computed unreliably anyway, ensuring transience of their errors
% may matter less than ensuring transience of errors in the (unreliable
% copy of the) sparse matrix and initial residual.  

%%%%%%%%%%%%%%%%%%%%%%%%%%%%%%%%%%%%%%%%%%%%%%%%%%%%%%%%%%%%%%%%%%%%%%
\subsubsection{Type system model}
%%%%%%%%%%%%%%%%%%%%%%%%%%%%%%%%%%%%%%%%%%%%%%%%%%%%%%%%%%%%%%%%%%%%%%

Consider implementing sparse matrix-vector multiply (the example of
Section \ref{SS:reliability-models:sandbox}) as the guest program in
the unreliable sandbox.  If the guest can be arbitrarily unreliable,
the sandbox has to do a lot of work to protect the host from things
like invalid instructions (due to errors in instructions) or
out-of-bounds array accesses (due to errors in index data).  The
sandbox could be much simpler if, for example, only the entries of the
sparse matrix and vectors, and the floating-point computations with
the matrix and vector values, are allowed to experience errors.  This
restriction would also make it easier for programmers to reason about
what happens in code running inside the sandbox, so they would not
need to write many redundant-looking checks that make code hard to
read and maintain.

\begin{comment}
%
% mfh 11 Jan 2011: Not a good example, because this works just fine
% in the "sandbox" model.
%
% mfh 23 Mar 2011: It works in the sandbox model if you decrease the
% scope of the "sandbox."
%
Here is another example: the accuracy of the sparse matrix-vector
products and preconditioner applications in an iterative method may
often be relaxed, as long as the basis vectors are kept orthogonal.
This is the principle behind ``inexact Krylov methods'' (see e.g.,
Simonici and Szyld \cite{simonici2003theory}), whose accuracy depends
on the orthogonality of the basis vectors, more so than in traditional
Krylov methods.  We might take this as inspiration to make the basis 
vectors reliable, and the sparse matrix-vector products and 
preconditioner applications unreliable.
\end{comment}

% , have more predictable performance, and
% perhaps even have predictable error bounds, 

This example suggests a finer-grained programming model, in which
developers can decide which data and computations they want to be
reliable or unreliable, and mix the two in their program.  For safety
and ease of use, the default behavior of all data and computations
should be as close to fail-stop reliability as possible.  (That is,
either the data and computations are reliable, or the program
terminates.)  Programmers may then relax reliability for certain data,
or certain phases of computation, or both.\footnote{Note that assuming
  a policy of default reliability and explicit unreliability does not
  contradict our characterization of bulk vs.\ high reliability.  It
  simply makes annotation easier.}  In the above example, fail-stop
default reliability ensures correctness of the sparse matrix indices
and the sparse matrix-vector multiply routine, so the routine will not
crash the entire program.  This programming model is more demanding
than the sandbox model, because it complicates the ways in which
reliable and unreliable computations and data may interact. 

We are currently exploring a special case of this model, in which
programmers can allocate ``unreliable memory'' by calling a special
version of C's \texttt{malloc} routine.  The operating system records
and reports to the application any detected but uncorrectible memory
faults in memory areas marked unreliable, but it does not kill the
process that allocated this memory, as many operating systems do for
ordinary memory allocations.  We believe this programming interface -
based approach will work for special cases of faults.  However, we
think the best way to generalize this reliability model for all kinds
of faults in different hardware components would be to encode
reliability in the type system of the programming language, much as
existing type systems encode the precision of floating-point values or
whether an object should be protected from simultaneous access by
multiple threads.  We do not require new programming language features
for the numerical methods proposed in this paper, but we think it
would make designing and implementing fault-tolerant algorithms much
easier.

\begin{comment}
Existing type systems in programming languages let
coders declare the precisions of floating-point objects, and define
type demotion or promotion rules when combining values of different
precisions.  Analogously, the type system could let coders declare
reliability of data explicitly as part of its type.  This in turn
would define which computations must be reliable, and how reliable and
unreliable data are allowed to mix in the same program.  For example,
addition of two reliable floating-point values should be computed
reliably, but addition of a reliable and a unreliable value may be
computed unreliably and result in a unreliable value (unless
explicitly promoted to a reliable value).  The resulting type analysis
may be done at compile time, run time, or both, as needed.  The
analysis need not be any more complicated than what existing compilers
already do when programmers mix different floating-point types.
Furthermore, the type system model may also be used to help make
errors look transient.  For example, arrays declared ``unreliable''
may be stored in unreliable main memory, but refreshed periodically
from a backing store.
\end{comment}

Encoding reliability in the type system is not a new idea.  Chen et
al.\ \cite{chen2005compiler} observe that different data in different
algorithms may need different levels of storage reliability, and that
reliability costs energy, space, performance, or some combination of
them.  They propose programmer annotations for declaring reliability
of subsets of multidimensional arrays.  For the simple case of nested
\texttt{for} loops over the arrays, they then use compiler analysis to
derive what parts of the arrays should be stored reliably.  Our
suggested ``reliability on demand'' feature is also a kind of
programmer annotation.  However, it applies to entire data structures
and computations, rather than subsets of arrays.  Chen et al.\ require
complicated compiler analysis of loops to derive the reliable regions
of arrays and generate separate reliable and unreliable code.  Our
annotations would depend only on simple type declarations and compiler
analysis, analogous to that already performed by compilers when
combining values of different floating-point precisions.  
% We also would not require complicated type deduction.

\begin{comment}
The above describes the type system reliability model in terms of the
programming language.  However, implementing it need not require
changes to existing programming languages.  It could be exposed to
programmers as a library instead, or by using object-oriented
programming language features like operator overloading.
\end{comment}

%%%%%%%%%%%%%%%%%%%%%%%%%%%%%%%%%%%%%%%%%%%%%%%%%%%%%%%%%%%%%%%%%%%%%%
\subsubsection{Reliable parallel decisions}
%%%%%%%%%%%%%%%%%%%%%%%%%%%%%%%%%%%%%%%%%%%%%%%%%%%%%%%%%%%%%%%%%%%%%%

Parallel computing introduces new ways in which soft faults can turn
into hard faults.
\begin{comment}
The above type system model lets programmers mix reliable and
unreliable data and computations in the same program.  This can
introduce new failure modes, depending on the allowed failure modes of
the unreliable parts.  
\end{comment}
For example, if the contents of messages between nodes of a
distributed-memory computer may become corrupted, then different nodes
may get different results in an all-reduce, even if each node computes
its part of the all-reduce reliably.  Many distributed-memory
implementations of iterative methods use the result of an all-reduce
in a predicate that tells the method when to stop iterating (for
example, when the residual norm is less than some tolerance).  The
predicate is computed redundantly on each node, with the expectation
that all nodes will get the same result.  If they do not -- for
example, if they have different values for the residual norm -- then
some nodes may stop iterating while others continue.  This can result
in deadlock or application failure, that is, it can turn a soft fault
into a hard fault.  We would prefer that \emph{parallel decisions}
like this one be reliable.

This is not a new problem; Blackford et al.\
\cite{blackford1996practical} discuss it in the less extreme context
of heterogeneous clusters, where different processors may have
different floating-point properties and thus may evaluate
floating-point comparisons differently.  They recommend in this case
that one processor compute the stopping criterion and broadcast the
Boolean result to all other processors.  This would only solve the
reliability problem for convergence tests if Boolean-valued messages
cannot be corrupted or lost.  In our case, it would be simpler, and
probably no more costly, to require the original all-reduce and the
predicate evaluation to be reliable and produce the same result on all
nodes.
\begin{comment}
Another approach would be to improve reliability of iterative method
coefficients (like the upper Hessenberg matrix in GMRES) and of the
computation of vector norms (even if the vector data itself on each
node is not reliable).  
\end{comment}

A different approach would be to observe that the stopping criterion
is a special case of distributed agreement on a Boolean value.  This
is an instance of the thoroughly studied Byzantine Generals Problem
(Lamport et al.\ \cite{lamport1982byzantine}), for which practical
solution algorithms exist (see e.g., Castro and Liskov
\cite{castro2002practical}).  This problem assumes that some of the
entities participating in distributed agreement may intentionally
attempt to deceive the others, which is an extreme but valid
generalization of corrupted data and arithmetic on some processors.
In practice, simple distributed agreement schemes should suffice.  For
example, an implementation could augment the all-reduce input for the
convergence test with a simple integer variable which each processor
would set to one if it has reached convergence.  Then all processors
would declare convergence if the sum of these integer values was
greater than some portion of the total processors being used.
Alternately, it may be simpler just to assume high reliability for all
distributed-memory transactions.  For example, practically speaking,
the cost of an all-reduce is dominated by latency (or even just the
fact that the message is transmitted off the node), so adding
reliability by computing redundantly or by adding error detection and
correction metadata to the all-reduce data package is almost free.

%% file: desired-properties.tex
\section{Desired properties of fault-tolerant iterations}\label{S:prop}
%%%%%%%%%%%%%%%%%%%%%%%%%%%%%%%%%%%%%%%%%%%%%%%%%%%%%%%%%%%%%%%%%%%%%%

Fault-tolerant iterative methods should have certain properties in
order to be both useful and feasible to implement.  In this section,
we describe a few desired properties, and explain which make sense to
implement.  Section \ref{SS:prop:conv} introduces two desired
convergence properties -- eventual convergence and gradual degradation
of convergence -- and argues for eventual convergence as the most
reasonable criterion.  Section \ref{SS:prop:impl} discusses properties
of implementations of these methods that will help them achieve good
performance, with minimal changes to existing solver algorithms and
implementations.  These criteria will help us narrow the space of
possible algorithms.

%%%%%%%%%%%%%%%%%%%%%%%%%%%%%%%%%%%%%%%%%%%%%%%%%%%%%%%%%%%%%%%%%%%%%%
\subsection{Convergence-related properties}\label{SS:prop:conv}
%%%%%%%%%%%%%%%%%%%%%%%%%%%%%%%%%%%%%%%%%%%%%%%%%%%%%%%%%%%%%%%%%%%%%%

We call what we see as the most important property \emph{eventual
  convergence}: If a comparable but not fault-tolerant method would
converge to the right answer in the case of no faults, the
fault-tolerant solver should either converge to the right answer in a
finite number of steps, or tell the caller that it did not.  The
fault-tolerant method may require more iterations or otherwise take
more time, and it might also have an upper bound on the number or
magnitude of faults it can tolerate.
\begin{comment}
% Ideally this
% property should hold in finite-precision arithmetic as well as in
% exact arithmetic.
\end{comment}
One iterative method that does \emph{not} have the eventual
convergence property is iterative refinement (an algorithm first
described by Wilkinson \cite{wilkinson1963rounding}).  Given
sufficiently large faults, only one fault in the residual vector need
happen at the ``last iteration'' for iterative refinement never to
compute the right answer.  Without eventual convergence, it would not
be worthwhile to relax hardware reliability, since all the effort at
previous iterations might be wasted by a single fault.  It is often
impossible to know when an fault will occur in a particular component,
so a reasonable method should allow them to occur at any time.  The
algorithm we present in Section \ref{S:alg} does have the eventual
convergence property.

\emph{Gradual degradation of convergence} as the number of faults
increases would also be desirable.  This might be much harder to
guarantee than eventual convergence.  For example, consider an
explicit Petrov-Galerkin projection method for solving the $n \times
n$ system $Ax=b$, that adds basis vectors to two different bases $V_k
= [v_1, \dots, v_k]$ and $W_k = [w_1, \dots, w_k]$.  Implementing a
method mathematically equivalent to GMRES, for instance, would require
$\Range{V_k}$ $=$ $\Span\{ r_0$, $A r_0$, $\dots$, $A^{k-1}$ $r_0 \}$
and $\Range{W_k} = A \Range{V_k}$.  If the matrix-vector products were
unreliable, we could still extend the basis in every iteration by
adding a random basis vector and orthogonalizing it against the
previous basis vectors, if the basis vectors are computed reliably.
In the worst case, this unreliable method would not converge until
$\Range{W_k}$ spans the entire space, that is, on iteration $n-1$.  In
fact, GMRES cannot promise better than this even in the case of no
faults.  It is possible to construct $n \times n$ problems for which
the residual in ordinary GMRES does not decrease until iteration
$n-1$, or for which the residual exhibits any desired nonincreasing
convergence curve \cite{greenbaum1996any}.  Some real-life linear
systems exhibit almost no convergence until some number of iterations,
after which they converge rapidly.  This suggests that eventual
convergence is a more reasonable goal than gradual degradation of
convergence.  We will show in the numerical experiments in Section
\ref{S:num-exp} that our FT-GMRES algorithm exhibits gradual
degradation of convergence in practice.  It may do so in theory also,
though we do not attempt in this paper to prove this.

%%%%%%%%%%%%%%%%%%%%%%%%%%%%%%%%%%%%%%%%%%%%%%%%%%%%%%%%%%%%%%%%%%%%%%
\subsection{Implementation-related properties}\label{SS:prop:impl}
%%%%%%%%%%%%%%%%%%%%%%%%%%%%%%%%%%%%%%%%%%%%%%%%%%%%%%%%%%%%%%%%%%%%%%

We have already discussed different models of application-controlled
reliability in Section \ref{S:reliability-models}.  Making all data
and arithmetic reliable would trivially result in a fault-tolerant
iterative method.  However, all of our models assume that reliability
has a cost, which is some combination of additional energy or storage
and reduced performance.  Thus, a fault-tolerant algorithm should aim
to store most of its data and spend most of its computations in
unreliable mode.  Second, fault-tolerant algorithms should not be too
much slower than corresponding less tolerant algorithms.  It is
reasonable to expect that the longer an application runs, the more
faults it will likely encounter.  More faults mean either slower
convergence, which compounds the problem, or even solver failure.  If
the fault-tolerant method is too slow, it may be faster just to run a
less tolerant method over and over using an ensemble approach until
the majority of answers agree.
\begin{comment}
An analogy would be asynchronous / chaotic
iterations: the algorithms themselves are asymptotically slower than
the best algorithms that do not tolerate message delays, so it is more
cost-effective (in current systems at least) to make messages reliable
and use faster algorithms.  
\end{comment}
Finally, fault-tolerant methods should reuse existing algorithms and
implementations as much as possible.  In particular, they should
accept existing preconditioner algorithms, and ideally even existing
implementations.  Preconditioners are often complicated and specific
to their application.  Our inner-outer iteration in Section
\ref{S:alg} can call existing iterative solvers and their
preconditioners as a ``black box,'' as long as they promise to
terminate within a fixed time.\footnote{Guaranteeing fixed-time
  termination when distributed-memory messages may be unreliable may
  require some modifications to existing sparse matrix-vector multiply
  and preconditioner implementations, but not to the mathematical
  algorithms.}

\begin{comment}
In practice, iterative methods can be more robust with respect to bad
inputs than their authors intended.  For example, I was recently
working with a physicist who was using CG to solve symmetric
indefinite systems.  His convergence plots showed spiky intermediate
behavior, which he correctly identified as ``deflating the negative
eigenvalues,'' but eventual convergence.  CG of course makes no
convergence guarantees in the case of an indefinite matrix, though
there is some theory to support its robustness: The CG tridiagonal
matrix may be ill-conditioned for one iteration, but never for two
consecutive iterations: see e.g.,
\cite{bank1993analysis,greenbaum1999solving}.
\end{comment}

%%% Local Variables: 
%%% mode: latex
%%% TeX-master: "paper"
%%% End: 

%% file: alg.tex
%%%%%%%%%%%%%%%%%%%%%%%%%%%%%%%%%%%%%%%%%%%%%%%%%%%%%%%%%%%%%%%%%%%%%%
\section{Fault-Tolerant GMRES}\label{S:alg}
%%%%%%%%%%%%%%%%%%%%%%%%%%%%%%%%%%%%%%%%%%%%%%%%%%%%%%%%%%%%%%%%%%%%%%

This section describes the Fault-Tolerant GMRES (FT-GMRES) algorithm,
a Krylov subspace method for iterative solution of large sparse linear
systems $Ax = b$.  FT-GMRES computes the correct solution $x$ even if
the system experiences uncorrected faults in both data and
arithmetic~\cite{hoemmen2011fault}.  It promises ``eventual
convergence'' in the sense of Section \ref{SS:prop:conv}: it will
always either converge to the right answer, or (in rare cases) stop
and report immediately to the caller if it cannot make progress.
FT-GMRES accomplishes this by dividing its computations into
\emph{reliable} and \emph{unreliable} phases, using the sandbox model
of reliability described in Section
\ref{SS:reliability-models:sandbox}.  Rather than rolling back any
faults that occur in unreliable phases, as a checkpoint / restart
approach would do, FT-GMRES ``rolls forward'' through any faults in
unreliable phases, and uses the reliable phases to drive convergence.
FT-GMRES can also exploit fault detection
%, including but not limited
% to the ECC memory fault notification discussed above, 
in order to correct corrupted data during unreliable phases.

%%%%%%%%%%%%%%%%%%%%%%%%%%%%%%%%%%%%%%%%%%%%%%%%%%%%%%%%%%%%%%%%%%%%%%
\subsection{FT-GMRES is based on Flexible GMRES}\label{SS:alg:FGMRES}
%%%%%%%%%%%%%%%%%%%%%%%%%%%%%%%%%%%%%%%%%%%%%%%%%%%%%%%%%%%%%%%%%%%%%%

\begin{algorithm}
\caption{(Right-preconditioned) GMRES}
\label{alg:GMRES}
\begin{algorithmic}[1]
\Input{Linear system $Ax=b$ and initial guess $x_0$}
\Output{Approximate solution $x_m$ for some $m \geq 0$}
\State{$r_0 := b - A x_0$}\Comment{Unpreconditioned initial residual vector}
\State{$\beta := \| r_0 \|_2$}
\State{$q_1 := r_0 / \beta$}
\For{$j = 1, 2, \dots$ until convergence}
  \State{Solve $q_j = M z_j$ for $z_j$}\Comment{Apply preconditioner}
  \State{$v_{j+1} := A z_j$}\Comment{Apply the matrix $A$}
  \For{$i = 1, 2, \dots, k$}\Comment{Orthogonalize}
    \State{$H(i,j) := q_i^* v_{j+1}$}
    \State{$v_{j+1} := v_{j+1} - q_i H(i,j)$}
  \EndFor
  \State{$H(j+1,j) := \| v_{j+1} \|_2$}
  \State{$q_{j+1} := v_{j+1} / H(j+1,j)$}\Comment{New basis vector}
  \State{$y_j := \argmin_y \| H(1:j+1,1:j) y - \beta e_1 \|_2$}
  \State{$x_j := x_0 + M^{-1} [q_1, q_2, \dots, q_j] y_j$}\Comment{Compute solution update}
\EndFor
\end{algorithmic}
\end{algorithm}

\begin{algorithm}
\caption{Flexible GMRES (FGMRES)}
\label{alg:FGMRES}
\begin{algorithmic}[1]
\Input{Linear system $Ax=b$ and initial guess $x_0$}
\Output{Approximate solution $x_m$ for some $m \geq 0$}
\State{$r_0 := b - A x_0$}\Comment{Unpreconditioned initial residual}
\State{$\beta := \| r_0 \|_2$}
\State{$q_1 := r_0 / \beta$}
\For{$j = 1, 2, \dots$ until convergence}
  \State{Solve $q_j = M_j z_j$}\label{alg:FGMRES:inner}\Comment{Apply current preconditioner}
  \State{$v_{j+1} := A z_j$}\Comment{Apply the matrix $A$}
  \For{$i = 1, 2, \dots, k$}\Comment{Orthogonalize}
    \State{$H(i,j) := q_i^* v_{j+1}$}
    \State{$v_{j+1} := v_{j+1} - q_i H(i,j)$}
  \EndFor
  \State{$H(j+1,j) := \| v_{j+1} \|_2$}
  \State{Update rank-revealing decomposition of $H(\IndexRange{1}{j}, \IndexRange{1}{j})$}
  \If{$H(j+1,j)$ is less than some tolerance}
    \If{$H(\IndexRange{1}{j},\IndexRange{1}{j})$ not full rank}
      \State{Did not converge; report error}\label{alg:FGMRES:rank-def}
    \Else
      \State{Return at end of this iteration}\Comment{Discovered invariant subspace}
    \EndIf
  \Else
    \State{$q_{j+1} := v_{j+1} / H(j+1,j)$}
  \EndIf
  \State{$y_j := \argmin_y \| H(\IndexRange{1}{j+1},\IndexRange{1}{j}) y - \beta e_1 \|_2$}
  \State{$x_j := x_0 + [z_1, z_2, \dots, z_j] y_j$}\Comment{Compute solution update}
\EndFor
\end{algorithmic}
\end{algorithm}

FT-GMRES is based on Flexible GMRES (FGMRES) \cite{saad1993flexible},
shown here as Algorithm \ref{alg:FGMRES}.  FGMRES extends the
Generalized Minimal Residual (GMRES) method of Saad and Schultz
\cite{saad1986gmres}, by ``flexibly'' allowing the preconditioner to
change in every iteration.  We show standard right-preconditioned
GMRES as Algorithm \ref{alg:GMRES} for comparison.  An important
motivation of flexible methods are ``inner-outer iterations,'' which
use an iterative method itself as the preconditioner.  In this case,
``solve $q_j := M_j z_j$'' (Line \ref{alg:FGMRES:inner} of Algorithm
\ref{alg:FGMRES}) means ``solve the linear system $A z_j = q_j$
approximately using a given iterative method.''  This \emph{inner
  solve} step preconditions the \emph{outer solve} (in this case the
FGMRES flexible iteration).  Changing right-hand sides and possibly
changing stopping criteria for each inner solve mean that if one could
express the ``inner solve operator'' as a matrix, it would be
different on each invocation.  This is why inner-outer iterations
require a flexible outer solver.

Flexible methods need not use an iterative method for the inner
solves.  The $M_j$ may be arbitrary functions from the range of $A$ to
the domain of $A$.  Most importantly, the preconditioners may change
significantly from one iteration to another; flexible methods do not
depend on the difference between successive preconditioners being
small.  This is the key observation behind FT-GMRES: flexible
iterations allow successive inner solves to differ arbitrarily, even
unboundedly.  This suggests modeling faulty inner solves as
``different preconditioners.''  Taking this suggestion leads to
FT-GMRES, which we present in the next section.

Flexible inner-outer iterations have the property that the dimension
of the Krylov subspace from which they choose the current approximate
solution grows at each outer iteration \cite{simonici2003flexible}, as
long as the break-down condition mentioned above does not occur.  This
ensures eventual convergence.  Corresponding restarted Krylov methods
lack this property; their convergence may stagnate.  Even though this
property of inner-outer iterations may not hold in the case of faulty
inner solves, the numerical experiments in Section \ref{S:num-exp}
show that inner-outer iterations offer better fault tolerance than
simply restarting.  Both restarting and inner-outer iterations
correspond naturally to the sandbox reliability model when the number
of iterations per restart cycle resp.\ inner solve is fixed.

There are flexible versions of other iterative methods besides GMRES,
such as CG \cite{golub1999inexact} and QMR \cite{szyld2001fqmr}, which
could also be used as the outer solver.  We chose FGMRES because it is
easy to implement, robust, and can handle nonsymmetric linear systems.
Experimenting with other flexible outer iterations is future work.

\subsubsection{Flexible GMRES' additional failure mode}

FGMRES has an additional failure mode beyond those of standard GMRES.
The quantity $H(j+1,j) = 0$ does not necessarily indicate convergence,
as it would in standard GMRES.  This is because $H(\IndexRange{1}{j},
\IndexRange{1}{j})$ is always nonsingular in GMRES if $j$ is the
smallest iteration index for which $H(j+1,j) = 0$, whereas
$H(\IndexRange{1}{j}, \IndexRange{1}{j})$ may not be nonsingular in
FGMRES in that case.  (This is Saad's Proposition 2.2
\cite{saad1993flexible}.)  This can happen even in exact arithmetic.
For example, it may occur due to unlucky choices of the
preconditioners: for example, $M_j^{-1} = A$ and $M_{j+1}^{-1} =
A^{-1}$.  In practice, this case is rare, even when inner solves are
subject to faults.  Furthermore, it can be detected inexpensively,
since there are algorithms for updating a rank-revealing decomposition
of an $m \times m$ matrix in $O(m^2)$ time (see e.g., Stewart
\cite{stewart1993updating}).  This is no more time than it takes to
update the QR factorization of the upper Hessenberg matrix at
iteration $m$.  The ability to detect this rank deficiency ensures
``trichotomy'' of FGMRES: it either
\begin{enumerate}
\item converges to the desired tolerance,
\item correctly detects an invariant subspace, with a clear indication
  ($H(j+1,j) = 0$ and $H(\IndexRange{1}{j}, \IndexRange{1}{j})$ is
  nonsingular), or
\item gives a clear indication of failure ($H(j+1,j) \neq 0$ and
  detected rank deficiency of $H(\IndexRange{1}{j},
  \IndexRange{1}{j})$).
\end{enumerate}
We base FT-GMRES' ``eventual convergence'' on this trichotomy
property.  In the following section, we will discuss recovery
strategies that FT-GMRES can use in case of the third condition above.

%%%%%%%%%%%%%%%%%%%%%%%%%%%%%%%%%%%%%%%%%%%%%%%%%%%%%%%%%%%%%%%%%%%%%%
\subsection{Fault-Tolerant GMRES}\label{SS:alg:FT-GMRES}
%%%%%%%%%%%%%%%%%%%%%%%%%%%%%%%%%%%%%%%%%%%%%%%%%%%%%%%%%%%%%%%%%%%%%%

FGMRES' acceptance of significantly different preconditioners at each
iteration suggests modeling solver faults as ``different
preconditioners.''  The least disruptive approach for existing solvers
is to use the inner-outer iteration approach.  The outer FGMRES
iteration wraps any existing solver, which is used as the inner
iteration.  Any solver works, even a sparse direct method (in which
case the inner ``iteration'' is not actually an iterative method), an
iterative method with any or no preconditioner, or a specialized
algorithm that exploits problem structure (such as an FFT or
hierarchical matrix factorization).  Existing preconditioners may also
be used without algorithmic modifications.  We call the resulting
inner-outer iteration \emph{Fault-Tolerant GMRES}.  It is shown here
as Algorithm \ref{alg:FT-GMRES}.  Inner-outer iterations with FGMRES
have been used as a kind of iterative refinement in mixed-precision
computation (see Buttari et al.\ \cite{buttari2006computations}), but
as far as we know, this is the first time it has been used for
reliability and robustness against possibly unbounded errors.

\begin{algorithm}
\caption{Fault-Tolerant GMRES (FT-GMRES)}
\label{alg:FT-GMRES}
\begin{algorithmic}[1]
\Input{Linear system $Ax=b$ and initial guess $x_0$}
\Output{Approximate solution $x_m$ for some $m \geq 0$}
\State{$r_0 := b - A x_0$}
\State{$\beta := \| r_0 \|_2$}
\State{$q_1 := r_0 / \beta$}\Comment{Unpreconditioned initial residual vector}
\For{$j = 1, 2, \dots$ until convergence}
  \State{Inner solve (unreliable) for $z_j$ in $q_j = A z_j$}\label{alg:FT-GMRES:inner}
  \State{$v_{j+1} := A z_j$}\Comment{Apply the matrix $A$}
  \For{$i = 1, 2, \dots, k$}\label{alg:FT-GMRES:orthog}\Comment{Orthogonalize}
    \State{$H(i,j) := q_i^* v_{j+1}$}
    \State{$v_{j+1} := v_{j+1} - q_i H(i,j)$}
  \EndFor
  \State{$H(j+1,j) := \| v_{j+1} \|_2$}
  \State{Update rank-revealing decomposition of $H(\IndexRange{1}{j}, \IndexRange{1}{j})$}
  \If{$H(j+1,j)$ is less than some tolerance}
    \If{$H(\IndexRange{1}{j},\IndexRange{1}{j})$ not full rank}
      \State{Try recovery strategies discussed in text}\label{alg:FT-GMRES:recovery}
    \Else
      \State{Return at end of this iteration}\Comment{Discovered invariant subspace}
    \EndIf
  \Else
    \State{$q_{j+1} := v_{j+1} / H(j+1,j)$}
  \EndIf
  \State{$y_j := \argmin_y \| H(\IndexRange{1}{j+1},\IndexRange{1}{j}) y - \beta e_1 \|_2$}
  \State{$x_j := x_0 + [z_1, z_2, \dots, z_j] y_j$}\Comment{Compute solution update}
\EndFor
\end{algorithmic}
\end{algorithm}

The only part of FT-GMRES allowed to run unreliably is Line
\ref{alg:FT-GMRES:inner}, which invokes the inner solver.  FT-GMRES
expects that inner solves do most of the work, so inner solves run in
the less expensive unreliable mode.  Inner solvers need only return
with a solution in finite time (see Section
\ref{SS:reliability-models:sandbox}).  That solution may be completely
wrong if errors occurrred.  Within the inner solves, the matrix $A$,
right-hand side $b$, and any other inner solver data may change
arbitrarily, and those changes need not even be transient.  However,
each outer iteration of FT-GMRES must run reliably, and requires a
correct version of the matrix $A$, right-hand side $b$, and additional
outer solve data (the same that FGMRES would use).

Since FT-GMRES expects only a small number of outer iterations,
interspersed by longer-running inner solves, we need not store two
copies (unreliable and reliable) of $A$ and $b$ in memory.  Instead,
we can save them to a reliable backing store, or even recompute them.
If the system provides fault detection capability, we can avoid
recovering or recomputing these data if no faults occurred, or even
selectively recover or recompute just the corrupted parts of the
critical data.  If the inner solve itself knows that no errors
occurred, it could also aggressively continue improving the solution
before returning to the outer iteration; we leave this option for
future work.

One practical point is that the outer iteration must scan the result
of each inner solve for invalid floating-point values (\texttt{NaN}
and \texttt{Inf}), and replace any with valid values.  The latter need
not be correct -- for example, they may be random numbers, or (better)
averages of their neighbors with respect to the matrix structure.
Many iterative methods perform this scan already for incomplete
factorization preconditioning, since there often is no way to know in
advance that the incomplete factors are nonsingular.

Line \ref{alg:FT-GMRES:recovery} of Algorithm \ref{alg:FT-GMRES}
covers the case where the outer iteration appears to have converged,
but the current upper Hessenberg matrix is rank deficient.  This can
happen in FGMRES as well, even with no faults.  There, it indicates an
unlucky combination of preconditioner applications.  In the case of
FT-GMRES, that unlucky combination may have occurred due to faults.
One of the following recovery strategies may be appropriate: 
\begin{enumerate}
\item retry the current iteration starting from Line
  \ref{alg:FT-GMRES:inner} inclusive;
\item retry the current iteration \emph{after} Line
  \ref{alg:FT-GMRES:inner}, but replace $z_j$ with a random
  vector (scaled appropriately according to best estimates of
  $\|A^{-1}\|$); or
\item stop and return $x_{j-1}$, the last good approximate solution.
\end{enumerate}
In parallel, all these strategies require agreement between
processors, and therefore global communication.  However, the
processors have to agree anyway whether to continue iterating based on
the convergence criterion, so no additional communication is needed.
In our numerical experiments discussed in Section \ref{S:num-exp}, we
found the rank-deficient upper Hessenberg case to be rare.

Another feature of the inner-outer iteration approach is that we can
reuse information from previous inner iterations, if we know somehow
that they were error-free.  For example, we could use a Krylov basis
recycling technique and restart, or simply keep the previous
iteration's data and continue without restarting (for an (F)GMRES
inner iteration).  Thus, the implementation can use whatever
information about errors is available, though it does not require this
information.

%%%%%%%%%%%%%%%%%%%%%%%%%%%%%%%%%%%%%%%%%%%%%%%%%%%%%%%%%%%%%%%%%%%%%%
\subsection{Inexact Krylov as an analysis tool}\label{SS:alg:inexact}
%%%%%%%%%%%%%%%%%%%%%%%%%%%%%%%%%%%%%%%%%%%%%%%%%%%%%%%%%%%%%%%%%%%%%%

\emph{Inexact Krylov methods} allow solving $Ax=b$ by using successive
approximations $A_k$ of $A$.  This makes them a generalization of
flexible methods, since the matrix, as well as the preconditioner, may
change in every iteration.  For overviews and development of
convergence theory, see Simonici and Szyld \cite{simonici2003theory}
and van den Eshof and Sleijpen \cite{eshof2004inexact}.  These methods
convergence when the error between the actual matrix $A$ and each
approximation $A_k$ respects a varying bound.  The bound starts small,
but grows inversely as a function of the current residual norm.
Inexact Krylov methods are motivated by applications where computing
$A$ itself is prohibitively expensive, but computing $w = A v$ for a
vector $v$ can be done approximately, and more effort in the
approximation results in less error.

Inexact Krylov methods cannot be used to provide tolerance against
arbitrary data and computational faults when applying the matrix $A$.
This is because they require an error bound which is usually not as
large as many possible bit flips.  (Bit flips may occur in exponent
bits as well as sign and significand bits.)  Furthermore, if a fault
in applying $A$ results in an error which is larger than the current
bound, inexact Krylov methods cannot promise convergence.
Nevertheless, inexact Krylov offers a framework for analyzing FGMRES
convergence.  If a reliability model lets us control and bound inner
solves' errors, we can use this framework.

Inexact Krylov methods also give insight into where to focus
reliability efforts.  For example, convergence of inexact GMRES
depends more on orthogonality of the basis vectors than convergence of
standard GMRES \cite{simonici2003theory}.  This suggests spending more
effort on basis vector reliability than on reliability of the matrix
and preconditioner.  Furthermore, the inexact Krylov framework
suggests that the matrix $A$ and preconditioner(s) should be applied
more reliably in initial iterations, if possible.  This coincides with
our informal experimental observation that perturbing the matrix $A$
affects convergence of iterative solvers more in earlier iterations
than in later iterations.

%%%%%%%%%%%%%%%%%%%%%%%%%%%%%%%%%%%%%%%%%%%%%%%%%%%%%%%%%%%%%%%%%%%%%%
%%%%%%%%%%%%%%%%%%%%%%%%%%%%%%%%%%%%%%%%%%%%%%%%%%%%%%%%%%%%%%%%%%%%%%

\begin{comment}
Finally, a
robust distributed-memory parallel implementation of the inner solve
should ensure that all processors agree on stopping criteria, in order
to prevent undefined behavior.  For a full discussion, see e.g.,
Blackford et al.\ \cite{blackford1996practical}, who consider this
problem in a different context.
\end{comment}

%%% Local Variables: 
%%% mode: latex
%%% TeX-master: "paper"
%%% End: 

%% file: programming-model.tex
%%%%%%%%%%%%%%%%%%%%%%%%%%%%%%%%%%%%%%%%%%%%%%%%%%%%%%%%%%%%%%%%%%%%%% 
\section{Programming model details}\label{S:model}
%%%%%%%%%%%%%%%%%%%%%%%%%%%%%%%%%%%%%%%%%%%%%%%%%%%%%%%%%%%%%%%%%%%%%%

When we presented the FT-GMRES algorithm in Section \ref{S:alg}, we
declared few assumptions about the reliability programming model.  The
algorithm needs few; the ``sandbox'' model (Section
\ref{SS:reliability-models:sandbox}) suffices for correctness, and
maps naturally to inner-outer iterations in general.  However,
existing computer systems require few modifications to offer a richer
model, which can also help us implement FT-GMRES more efficiently.  In
this section, we describe a programming model that is both suited for
FT-GMRES, and is reasonable for systems architects to implement.  This
model promises reliable data and computations within the specified
time and space bounds, and provides best-effort fault detection
outside those bounds.  It includes schemes for efficient local
recovery of possibly corrupted data.  We were able to implement a
representative subset of this model for our performance prototype
(Section \ref{S:perf-exp}) with reasonable effort.

% Computing reliably with unreliable hardware requires a programming
% model that delimits faults in time and space.  

In Section \ref{SS:model:which}, we show how the data in FT-GMRES and
analogous methods naturally separate into categories based on required
reliability, and the amount of time and memory it consumes.  Section
\ref{SS:model:detection} explains why we assume only best-effort fault
detection, though better fault detection guarantees could improve
performance.  In Section \ref{SS:model:recovery}, we describe how the
FT-GMRES algorithm itself, plus best-effort fault detection, lead to a
two-stage recovery scheme for corrupted data.  This scheme makes
approximate repairs in inner iterations, and performs preferably local
full recovery of corrupted data in outer iterations.  

Throughout this section, we refer to the subset of the model we
implemented for the performance experiments discussed in Section
\ref{S:perf-exp}.  This system currently only handles detected but
uncorrectible DRAM memory faults, not other kinds of faults such as
incorrect arithmetic or corrupted MPI messages.  This restriction was
convenient for developing a prototype in reasonable time.  We argue in
Section \ref{SS:model:restriction} that a system that considers only
memory faults nevertheless includes the right programming model
elements for developing algorithms that can handle all kinds of
faults.  Finally, Section \ref{SS:model:general} proposes that our
model is sufficiently general that it could work for other numerical
methods based on subspace search and fixed-point iteration.

%%%%%%%%%%%%%%%%%%%%%%%%%%%%%%%%%%%%%%%%%%%%%%%%%%%%%%%%%%%%%%%%%%%%%%
\subsection{Which data reliable, when}\label{SS:model:which}
%%%%%%%%%%%%%%%%%%%%%%%%%%%%%%%%%%%%%%%%%%%%%%%%%%%%%%%%%%%%%%%%%%%%%%

In this section, we explain which data in our fault-tolerant iterative
method we allow to experience faults, and when in the algorithm we
allow those faults to occur.  In particular, we allow faults in all
``large'' data and computations in the inner iterations only.
``Large'' data includes sparse matrices, preconditioners, and vectors,
but does not include the small projected linear system or
least-squares problem used to compute the solution update
coefficients, nor does it include code or control data such as loop
indices.  We also explain which of the large data require occasional
recovery to their original uncorrupted state.  Finally, we argue that
this programming model could apply to other Krylov subspace methods,
and to subspace search and fixed-point iterations in general.

\subsubsection{``Large'' and ``small'' data}

Krylov methods for solving linear systems spend most of their memory
and time computing with two kinds of objects: ``large'' dense vectors,
and linear operators (functions from a vector to a vector) of the same
dimension(s).  The latter include sparse matrices (where the function
is sparse matrix-vector multiplication), linear operators implemented
as a subroutine (e.g., by discretizing and solving a partial
differential equation) rather than as a sparse matrix, and
preconditioners (if applicable).  

Krylov methods project a larger linear system onto a smaller linear
system or least-squares problem which is inexpensive to solve using
either dense factorizations, or an equivalent small number of scalar
computations.  This gives us a subjective but practical definition of
``large'' data: using a Krylov method to solve a linear system of that
size saves time, memory, or both, relative to a dense factorization.
Krylov methods also include ``small'' data: scalars or small dense
matrices and vectors which represent the projected linear system or
least-squares problem.  The projected problem is used to solve for the
coefficients of the solution update.  The projected problem requires
little memory or solution time computed with the large vectors and
operators.  Its small size makes it sensitive to corrupted data or
computations, yet the resulting solution update coefficients have a
large effect on the accuracy of the computed solution vector.  Thus,
we require that the projected problem be stored and computed reliably,
and confine any unreliable data or computation to the large vectors
and operators.

\subsubsection{Both operators and vectors may be unreliable}

The large matrix $A$ and preconditioner(s) typically take up much more
memory than a single vector or corresponding size.  Also, applying the
matrix or a preconditioner to a vector takes more time than computing
a single vector operation (such as a norm, inner product, or vector
sum).  However, the balance of time and memory between operator
applications and vector operations varies between Krylov methods.  For
example, our inner solver uses GMRES (the Generalized Minimal Residual
method \cite{saad1986gmres}), which may spend more of its time in
vector operations, depending on the restart length.  Thus, we allow
vectors as well as operators in the inner solver to be unreliable,
since otherwise the solver might require too much unreliable data and
computation.  The goal is for a fault-tolerant solver to spend most of
its time and memory in unreliable mode.

%%%%%%%%%%%%%%%%%%%%%%%%%%%%%%%%%%%%%%%%%%%%%%%%%%%%%%%%%%%%%%%%%%%%%%
\subsection{Best-effort fault detection}\label{SS:model:detection}
%%%%%%%%%%%%%%%%%%%%%%%%%%%%%%%%%%%%%%%%%%%%%%%%%%%%%%%%%%%%%%%%%%%%%%

We pessimistically assume best-effort fault detection.  This means
that a significant fraction of faults might evade detection.  We
assume this in part due to technical limitations of our software
prototype.  Currently, it can only detect injected faults by
simulating an ECC memory ``patrol scrubber'' in software, using a
separate, asynchronously executing thread.  Injected faults
encountered by user code's actual memory operations are not detected.
This is because the current version of Linux, on which our software
prototype depends, kills any user process whose memory operations
encounter an uncorrectable fault.  (It need not do this for faults
detected by an actual patrol scrubber.)  Changing this behavior would
require a custom Linux patch, which in turn would prohibit us from
running tests on computers we do not administer.  Many other operating
systems have this property.

Despite this technical limitation of our prototype, we believe that
our pessimistic assumption is reasonable in production systems.  For
example, most current systems offer no hardware detection of
arithmetic faults.  Without expensive hardware replication, the best a
system could do is insert occasional test instructions into the
instruction stream.  This would be mostly likely to detect ``sticky''
arithmetic faults, but not transient ones.  Detection might also be
asynchronous, so that faults contaminate other computations
irreversibly before the system detects and reports them to the
application.  For instance, a hardware ECC memory patrol scrubber
might discover uncorrectable corruption in a sparse matrix entry while
the Krylov method is doing something else.  Future systems may not
necessarily promise anything about the delivery time of the resulting
error report.  Finally, fault detection does cost energy and / or
performance.  Our algorithm does not require infallible fault
detection for correctness, so we are willing to relax this, as long as
system architects can meet our reliability demands.

%%%%%%%%%%%%%%%%%%%%%%%%%%%%%%%%%%%%%%%%%%%%%%%%%%%%%%%%%%%%%%%%%%%%%%
\subsection{Repair of corrupted data}\label{SS:model:recovery}
%%%%%%%%%%%%%%%%%%%%%%%%%%%%%%%%%%%%%%%%%%%%%%%%%%%%%%%%%%%%%%%%%%%%%%

\subsubsection{How Krylov methods use operators}

Krylov methods for solving linear systems use ``large'' linear
operators in two different ways.  The first way is iterative: the
method repeatedly applies the operator(s) to a vector, in order to
build up one or more search spaces.  The theory of \emph{inexact
  Krylov methods} says that the operators need not be applied exactly
at all iterations in order for the method to converge.  We take this
as inspiration for allowing these operator applications to vary
arbitrary, due to unreliability.  The second way is for computing the
residual vector of the current approximate solution explicitly.  The
residual vector may then be used to verify the approximate solution,
restart the iteration, or improve the solution in an outer iteration.
These uses of the residual vector require an exact computation, not an
approximation, using just the matrix $A$ and right-hand side $b$.
Techniques like iterative refinement even require computing the
residual vector in higher precision, in order for certain convergence
results to hold.

Constructing the operators always happens outside of the Krylov
method.  Construction may be a complicated operation consuming a
significant part of the application's total run time, and many more
lines of code than the linear solver.  (Consider a structural dynamics
application using the finite element method, for example.)  It is
usually a \emph{nonlocal} operation as well: that is, it requires
communication when running in parallel in a distributed-memory
environment.  (For example, in the finite element method, assembling
elements with mesh points owned by different processors requires
summing contributions from the involved processors.)  However, the
operators usually do not change or need expensive reconstruction
during the Krylov method.  

Note that solving linear systems with a Krylov method often requires
an effective preconditioner.  Preconditioners can be time-consuming to
compute, and this computation often requires global communication.
Algebraic multigrid is a good example.  The only difference between a
preconditioner and the matrix $A$ is that for GMRES variants, a left
preconditioner is not needed in order to compute a residual vector.
Many GMRES users prefer a right preconditioner anyway, since it
ensures that the projected least-squares problem's solution has the
same residual norm as the approximate solution, in exact arithmetic.
GMRES requires the right preconditioner in order to compute the
approximate solution or current residual vector; thus, computing these
vectors reliably requires applying the right preconditioner reliably.

\subsubsection{What this says about operator recovery}

The above two paragraphs say that (a) a fault-tolerant Krylov method
must be able to apply operators both reliably and unreliably, and (b)
constructing an operator is expensive and nonlocal.  This suggests
that a fault-tolerant Krylov method must be able to recover the
original operators reliably, and that this recovery should not require
recomputing the affected operator.  We suggest implementing this using
a \emph{local checkpointing} scheme, which saves matrix and
preconditioner data that may experience memory faults to reliable
backing storage.  We already assume that FT-GMRES marks this data as
unreliable, and that FT-GMRES notifies the system on entry to each
inner iteration that faults are allowed.  The checkpointing mechanism
need only pay attention to these notifications to decide when and what
to checkpoint.  The reliable backing store should be fast,
nonvolatile, and local to each node.  Recent projections for
exascale-class systems predict much heavier use of node-local scratch
storage (see e.g., \cite[Section 5.6.3.1]{kogge:exa}).  We expect,
therefore, that future supercomputers will include node-local
solid-state drives, meant for scratch storage or as a cache for input
/ output operations.  We currently lack access to such hardware, so as
a proxy, we implemented for this paper a reliable backing store using
ordinary DRAM memory in which we do not allow detected but
uncorrectible ECC faults (see Section \ref{S:perf-exp}).

\subsubsection{Krylov basis vectors}

We did not mention the Krylov basis vectors computed by the inner
iteration in the paragraphs above.  These vectors result from applying
a possibly corrupted matrix or preconditioner; they are ``corrupted by
construction.''  Thus, it does not make sense to save or restore them.
Vectors computed by the outer iteration should be completely reliable,
however.  Corruption of Krylov basis vectors in the outer iteration
may result in an incorrect solution.  

\subsubsection{Local and approximate recovery}

\emph{Local} recovery is important.  Faults like bit flips in memory
and incorrect arithmetic are local to the node (or even to the CPU).
Recomputation of an operator typically involves global communication,
whose pattern of dependencies typically make it a heavyweight global
synchronization point.  As supercomputers grow towards exascale, the
increasing cost of communication makes favoring local operations more
and more attractive.  The checkpointing scheme mentioned above offers
an exact local recovery method.  If the system offers reliable
detection of data corruption, including fault locations,
\emph{approximate} local repair is possible.  As we explain in Section
\ref{S:perf-exp}, existing ECC memory hardware provides this
information to the operating system upon encountering an uncorrectible
error.  The application can then define a handler that repairs the
fault.  For example, a corrupted sparse matrix entry can be ``smoothed
out'' by replacing it with the average of its neighbors.  Simple
handlers cost little more than the system interrupt caused by the
fault itself.

Approximate repair is an inexpensive option for inner iterations.
However, outer iterations require exact recovery of operators.  Since
corrupted data locations are unknown in advance, restoring the
operators requires either full local checkpointing, or global
recomputation.  This suggests a two-fold recovery strategy.  Start
each inner iteration with the correct sparse matrix and
preconditioner(s), but allow data corruption to occur.  If possible,
try to fix corrupted values during the inner iteration, but do so only
locally, and as quickly as possible, even if that means the values are
only recovered approximately.  Quick fixes minimize the idle time of
other processors which did not experience data corruption.  Local
fixing avoids communication overhead in the performance-critical inner
iteration.  At the end of the inner iteration, refresh the correct
values in the sparse matrix, even if no faults were detected there.
This ensures correctness even if undetected faults occurred.  Perform
the outer iteration, and continue.  Since we expect outer iterations
to occur infrequently, we can afford to spend more there on recovery
than in inner iterations.

%%%%%%%%%%%%%%%%%%%%%%%%%%%%%%%%%%%%%%%%%%%%%%%%%%%%%%%%%%%%%%%%%%%%%%
\subsection{Summary of model}\label{SS:model:summary}
%%%%%%%%%%%%%%%%%%%%%%%%%%%%%%%%%%%%%%%%%%%%%%%%%%%%%%%%%%%%%%%%%%%%%%

The above discussion implies three tiers of data and computation:
\begin{enumerate}
\item Always reliable data, such as the projected linear system, code,
  and control data (e.g., loop indices).
\item Data which may be unreliable in inner phases, must be reliable
  in outer phases, and which the outer iteration must be able to
  refresh to correct values.  Examples: the matrix, preconditioner(s),
  and right-hand side of the linear system $Ax=b$ to solve.
\item Unreliable data which does not require saving or restoring, such
  as the Krylov basis vectors in inner iterations.
\end{enumerate}

\begin{comment}
The sparse matrix and preconditioner may consume most of the memory
for some Krylov methods.  Indeed, we expect this ratio of matrix and
preconditioner data to vector data to increase in the future.  Since
total memory per processor will likely decrease over time, users may
favor methods that use fewer vectors.  
\end{comment}

%%%%%%%%%%%%%%%%%%%%%%%%%%%%%%%%%%%%%%%%%%%%%%%%%%%%%%%%%%%%%%%%%%%%%%
\subsection{Memory faults are sufficient}\label{SS:model:restriction}
%%%%%%%%%%%%%%%%%%%%%%%%%%%%%%%%%%%%%%%%%%%%%%%%%%%%%%%%%%%%%%%%%%%%%%

The performance prototype of FT-GMRES we describe in this work was
designed to handle faults in DRAM memory.  Computer hardware may also
experience corrupted caches or registers, arithmetic computations, or
messages between processors.  (This paper only considers faults that
result in corrupted data; other fault-tolerance techniques apply to
events like dropped messages or crashed nodes.)  Nevertheless, we
think that that the above programming model, and our fault-tolerant
inner-outer iteration approach, apply more generally to all kinds of
faults.

%
% Sandbox model combined with type system model covers all kinds of faults.
%

\begin{comment}
Our programming model restricts faults in two different ways.  First,
they are only allowed to occur within a defined scope.  This is the
inner solver in our case.  Second, they are only allowed to occur to
certain data.
\end{comment}

%% Recovering from floating-point arithmetic faults may require global
%% backtracking.

Floating-point arithmetic faults differ from memory faults, in that
there is no storage location to recover to an original value.  Thus,
bounding them in space is impossible.  Local fault recovery doesn't
make sense, because there is no storage location to recover.  However,
bounding them in time \emph{is} possible; one can use any of various
hardware or software approaches (e.g., triple modular redundancy) to
do so.  Furthermore, a reliable outer iteration can correct the
\emph{effects} of arithmetic faults in inner iterations, using an
algorithmic approach.  Thus, our solver could be easily made tolerant
of floating-point arithmetic faults as well.

The possibility of corrupted distributed-memory messages would violate
the principle our model assumes, namely that faults are local.
However, corrupted messages can be changed from a global to a local
issue by using error-correcting codes.  Message-passing hardware often
does this anyway.  Such codes enable the receiver of a corrupted
message to recover its original contents without communication.  Since
the latency of sending messages over a network is slow anyway compared
with computation, it is worthwhile paying the computational and
message bandwidth cost of an error correction scheme.  Furthermore, we
can model some kinds of corrupted messages (for example, when
computing a distributed sparse matrix-vector multiply) as transient
corruption of the operator.  Iterative methods do require that
stopping criteria (which are global Boolean decisions) be computed
reliably.  See \cite{blackford1996practical,blackford1997practical}
for a discussion of this issue in the context of heterogeneous compute
nodes.  In practice, making stopping criteria robust has little
performance penalty.

\begin{comment}
Dropped messages and node crashes are a different case (``hard''
faults).  For a known communication pattern (e.g., not \verb!ANY_TAG!
for the source), it's easy to retransmit, so you can assume messages
won't get dropped.

Corruption may happen in the sender's memory buffers, before the
message even touches the network.

%%%%%%%%%%%%%%%%%%%%%%%%%%%%%%%%%%%%%%%%%%%%%%%%%%%%%%%%%%%%%%%%%%%%%%
\subsubsection{(Floating-point) arithmetic faults}
%%%%%%%%%%%%%%%%%%%%%%%%%%%%%%%%%%%%%%%%%%%%%%%%%%%%%%%%%%%%%%%%%%%%%%

Consequences may propagate globally (e.g., reduction result).

Detection can therefore only be region- / time phase - based.  This is
just like the worst case of memory fault detection.

Local recovery doesn't make sense, due to propagation, but
``smoothing'' / global techniques for after-the-fact correction make
more sense.  That's the point of outer iteration -- it is fancy
smoothing.

Arithmetic faults force a global, algorithmic approach.  It's easier
to treat storage errors with a local checkpoint / restart - like
technique.
\end{comment}

%%%%%%%%%%%%%%%%%%%%%%%%%%%%%%%%%%%%%%%%%%%%%%%%%%%%%%%%%%%%%%%%%%%%%%
\subsection{Our model applies to other numerical methods}
\label{SS:model:general}
%%%%%%%%%%%%%%%%%%%%%%%%%%%%%%%%%%%%%%%%%%%%%%%%%%%%%%%%%%%%%%%%%%%%%%

%
% FIXME (mfh 02 Dec 2011) Perhaps this discussion should move
% elsewhere, to a more algorithm-oriented section.
%

Other numerical algorithms besides Krylov methods involve inner-outer
iterations based on repeatedly applying operators to vectors.
Newton's method and its variants for solving nonlinear equations are
one example.  In this case, the repeated linear solves form the inner
iterations.  
%
% FIXME (mfh 02 Dec 2011) Need reference for "practical
% implementations of Newton's method."
%
Practical implementations of Newton's method typically expect some
inner iterations to go awry, and ensure eventual convergence at the
outer level using trust region techniques.  In this case, the
operators and vectors in the linear solves may experience occasional
data or arithmetic faults.  The outer solves' residual, line search,
and trust region computations must be reliable.

Fixed-point iterations such as Picard iteration, so-called
``stationary iterative methods'' like Schwarz domain decomposition, or
even iterative refinement are other examples of inner-outer iterations
based on repeated applications of operators to vectors.  Depending on
the algorithm, these may or may not have guarantees of eventual
convergence in the presence of occasional faults.  Nevertheless, in
practice, the algorithms may still converge despite faults, so it
would be worthwhile exploring adding the fault model to them.

%%% Local Variables: 
%%% mode: latex
%%% TeX-master: "paper"
%%% End: 

%% file: numerical-experiments.tex
%%%%%%%%%%%%%%%%%%%%%%%%%%%%%%%%%%%%%%%%%%%%%%%%%%%%%%%%%%%%%%%%%%%%%%
\section{Numerical experiments}\label{S:num-exp}
%%%%%%%%%%%%%%%%%%%%%%%%%%%%%%%%%%%%%%%%%%%%%%%%%%%%%%%%%%%%%%%%%%%%%%

We began our experiments by prototyping solvers and a fault injection
framework in MATLAB.\footnote{MATLAB\textregistered is a registered
  trademark of The MathWorks, Inc.  We used MATLAB version 7.6.0.324
  (R2008a).}  We used these to compare the convergence of FT-GMRES,
restarted GMRES, and nonrestarted GMRES, for various fault rates in
the inner solves' sparse matrix-vector multiplies (SpMVs).  For these
experiments, we allowed only SpMV operations to experience faults, and
did not apply preconditioning.  Our performance prototype experiments
described in Section \ref{S:perf-exp} include preconditioning, and
allow faults anywhere in the inner solves.

Our initial experiments show that FT-GMRES can often converge even
when the majority of the inner solves' SpMVs suffer faults.  The other
methods tested either did not converge, or converged much more slowly
than FT-GMRES, when some of their SpMVs were faulty.  Furthermore,
FT-GMRES' convergence exhibits the desired gradual degradation
behavior as the fault rate increases.  Section
\ref{SS:num-exp:framework} describes our framework for numerical
experiments, and the test problems and actual experiments we tried.
We present results in Section \ref{SS:num-exp:results}.

%%%%%%%%%%%%%%%%%%%%%%%%%%%%%%%%%%%%%%%%%%%%%%%%%%%%%%%%%%%%%%%%%%%%%%
\subsection{Experimental framework}\label{SS:num-exp:framework}
%%%%%%%%%%%%%%%%%%%%%%%%%%%%%%%%%%%%%%%%%%%%%%%%%%%%%%%%%%%%%%%%%%%%%%

Our MATLAB prototype can inject faults either in the result of an
SpMV, or an entire inner solve (for FT-GMRES).  It decides
deterministically whether to inject a fault, by using a repeating
infinite sequence of Boolean values that we specify.  Each ``possibly
faulty'' operation reads the current Boolean value from the sequence,
and if it is true, we add 1 to the first entry of the result of the
operation (imitating \cite{howle2010soft:Copper}).  For example, when
running FT-GMRES with faulty SpMV operations, if the sequence is 0, 0,
1, then every third SpMV operation in the inner solve is faulty.
Deterministic faults make it easy to reproduce experimental results.
They also let us control which SpMV operations fail.  (This is
important because the theory of inexact Krylov methods (see Section
\ref{SS:alg:inexact}) suggests that inaccurate matrix-vector products
or preconditioner applications in the first few iterations matter more
than in later iterations.  We plan to explore this more in future
work.)

\begin{comment}
For simplicity, we always inject a fault by adding one to the first
entry of the output vector of the faulty operation.  (Our framework
can also replace the entire output vector with random data, to
simulate an aborted inner solve with serious memory corruption.)
\end{comment}

Our MATLAB versions of GMRES and FT-GMRES do extra work for
robustness.  After invoking a possibly unreliable operation (either an
SpMV or an inner solve), they scan the output vector for invalid
floating-point values (\texttt{Inf} or \texttt{NaN}), and replace
those with random data.  Also, after orthogonalization, they check
whether the norm of the resulting orthogonalized vector is an invalid
floating-point value.  If it is, they replace it with random data and
reorthogonalize.\footnote{Randomization improves robustness in
  practice, but makes reproducing experiments more difficult.  We used
  MATLAB's default Mersenne Twister pseudorandom number generator,
  with the default seed.}  Finally, we found that FT-GMRES converges
faster if the first inner solve is successful.  We implemented extra
reliability for the first inner solve in a realistic way as follows.
If the first inner solve did not reduce the residual norm at all, we
try it once more.  If that still did not reduce the residual norm, we
replace the result of the first inner solve with the identity operator
and continue.  We include this only for the first outer iteration of
FT-GMRES.  In practice, our experiments rarely needed to retry the
first inner solve.

We performed three sets of numerical experiments.  First, for a given
linear system and fault sequence, we compared the convergence of (a)
FT-GMRES, with $s - k + 1$ iterations per inner solve at outer
iteration $k$, for a total of $t$ outer iterations ($k = 1, \dots,
t$); (b) restarted GMRES, with $s$ iterations per restart cycle and
$t$ restart cycles; and (c) GMRES without restarting, $s \cdot t$
iterations.  Decreasing the number of iterations per inner solve in
FT-GMRES makes comparing an inner-outer iteration and a restarted
method fair, by ensuring that both methods store the same number of
left Krylov basis vectors \cite{saad1993flexible}.  We include
nonrestarted GMRES just to show its lack of robustness in the presence
of faults.  For this set, we fixed $s = 50$, to simulate the
fixed-time requirement for inner solves.  We set $t = 10$ so that $s
\cdot t$ nonrestarted GMRES iterations would complete in a reasonable
time.  Second, we tested only FT-GMRES with the same linear system,
but with different fault rates.  This set will show the desired
gradual degradation of FT-GMRES's convergence with respect to the
fault rate.  Here, we set $s-k+1$ iterations per inner solve with
$s=50$ as before, but performed more outer iterations ($t=20$), since
we did not have to run $s \cdot t$ iterations of nonrestarted GMRES.
In the third set, we tested FT-GMRES for many outer iterations $t=300$
and a fixed number $s=50$ of iterations per inner solve, and varied
the outer solves' convergence tolerance and the fault rate.  This will
show that computational cost does not increase much as the fault rate
increases.

\begin{table}
\begin{center}
\begin{tabular}{l|r|r|r}
Name           & \# rows & \# nz   & $\kappa(A)$ \\ \hline
Diagonal       & 10,000  & 10,000  & \verb!1.00e+10! \\
Szczerba/      & 20,896  & 191,368 & \verb!4.85e+09! \\
Ill\_Stokes    &         &         &                 \\
Sandia/        & 25,187  & 193,216 & \verb!1.99e+14! \\
mult\_dcop\_03 &         &         &                 \\
\end{tabular}
\end{center}
\caption{\label{tbl:matrices} Test problems for FT-GMRES numerical
  experiments.  The ``name'' (except for ``Diagonal'') comes from the
  University of Florida Sparse Matrix Collection.  ``\# rows'' gives 
  the number of rows (and columns) in the matrix, ``\# nz'' the number 
  of stored sparse matrix entries, and ``$\kappa(A)$'' an (estimate
  of, via MATLAB's \texttt{condest}) the matrix's condition number to 
  3 significant figures.}
\end{table}

We tested three types of matrices in our experiments: diagonal with
positive entries with base-10 logarithmic spacing from 1 to
$10^{-10}$, nonsymmetric matrices from discretizations of partial
differential equations (PDEs), and nonsymmetric circuit simulation
matrices.  Our matrices from the latter two categories come from the
University of Florida Sparse Matrix Collection (UFSMC) \cite{UFSMC}.
\begin{comment}
Diagonal matrices with positive entries serve as a proxy for symmetric
positive definite matrices; these are easy for GMRES to solve, so they
help isolate the effects of faults from the difficulty of the problem
itself.  PDE discretizations are popular test problems for new
iterative methods.  Circuit simulation problems often have a
fundamentally different structure and convergence behavior than
problems from the other two categories.  
\end{comment}
Table \ref{tbl:matrices} names and describes the test problems.
``Diagonal'' is a diagonal matrix, Ill\_Stokes comes from a
discretization of Stokes' equation, and mult\_dcop\_03 comes from a
circuit simulation.  Each UFSMC matrix includes a sample right-hand
side from its application.  
%
% For comparison, we first calculated the
% exact solution using a direct method.  
%
For ``Diagonal,'' we chose the exact solution $x$ as a vector of ones,
and computed the right-hand side $b$ via $b = A \cdot x$.

%%%%%%%%%%%%%%%%%%%%%%%%%%%%%%%%%%%%%%%%%%%%%%%%%%%%%%%%%%%%%%%%%%%%%%
\subsection{Results}\label{SS:num-exp:results}
%%%%%%%%%%%%%%%%%%%%%%%%%%%%%%%%%%%%%%%%%%%%%%%%%%%%%%%%%%%%%%%%%%%%%%

\begin{figure}
\begin{center}
\includegraphics[scale=0.55]{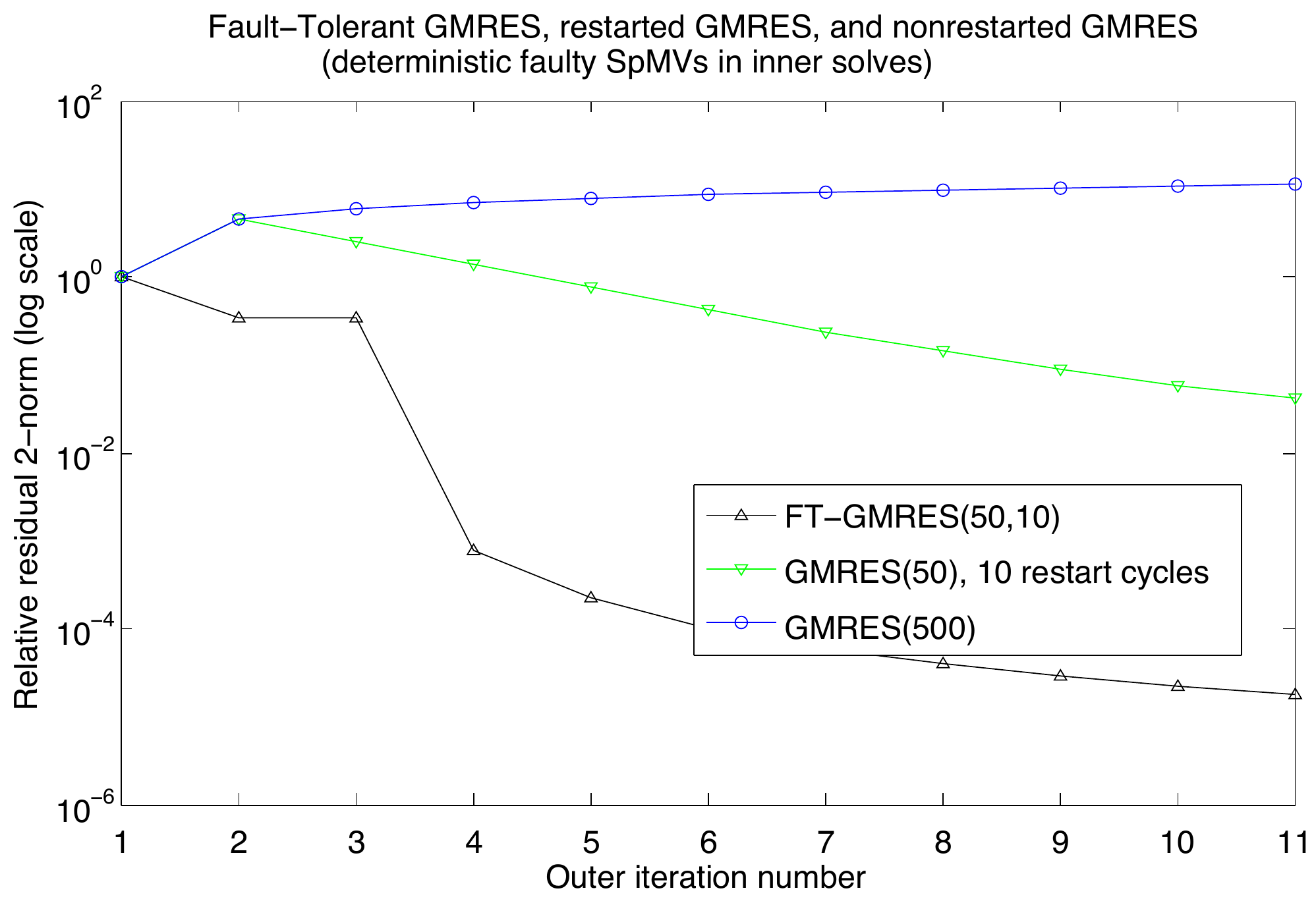}
\end{center}
\caption{\label{fig:diag10k:restarted} FT-GMRES vs.\ GMRES on Diagonal.}
\end{figure}

\begin{figure}
\begin{center}
\includegraphics[scale=0.5]{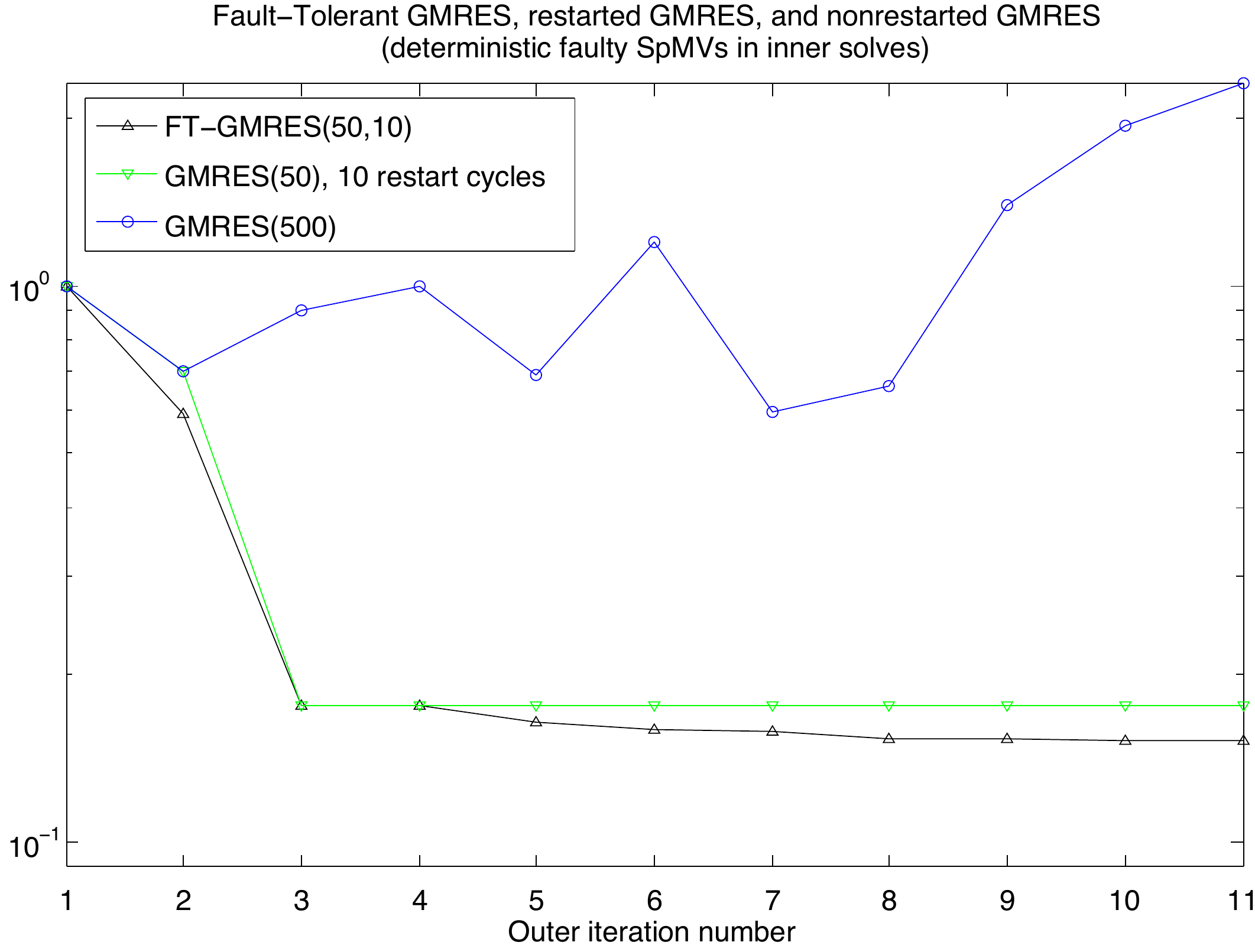}
\end{center}
\caption{\label{fig:IllStokes:restarted} FT-GMRES vs.\ GMRES
  on Ill\_Stokes.}
\end{figure}

\begin{figure}
\begin{center}
\includegraphics[scale=0.5]{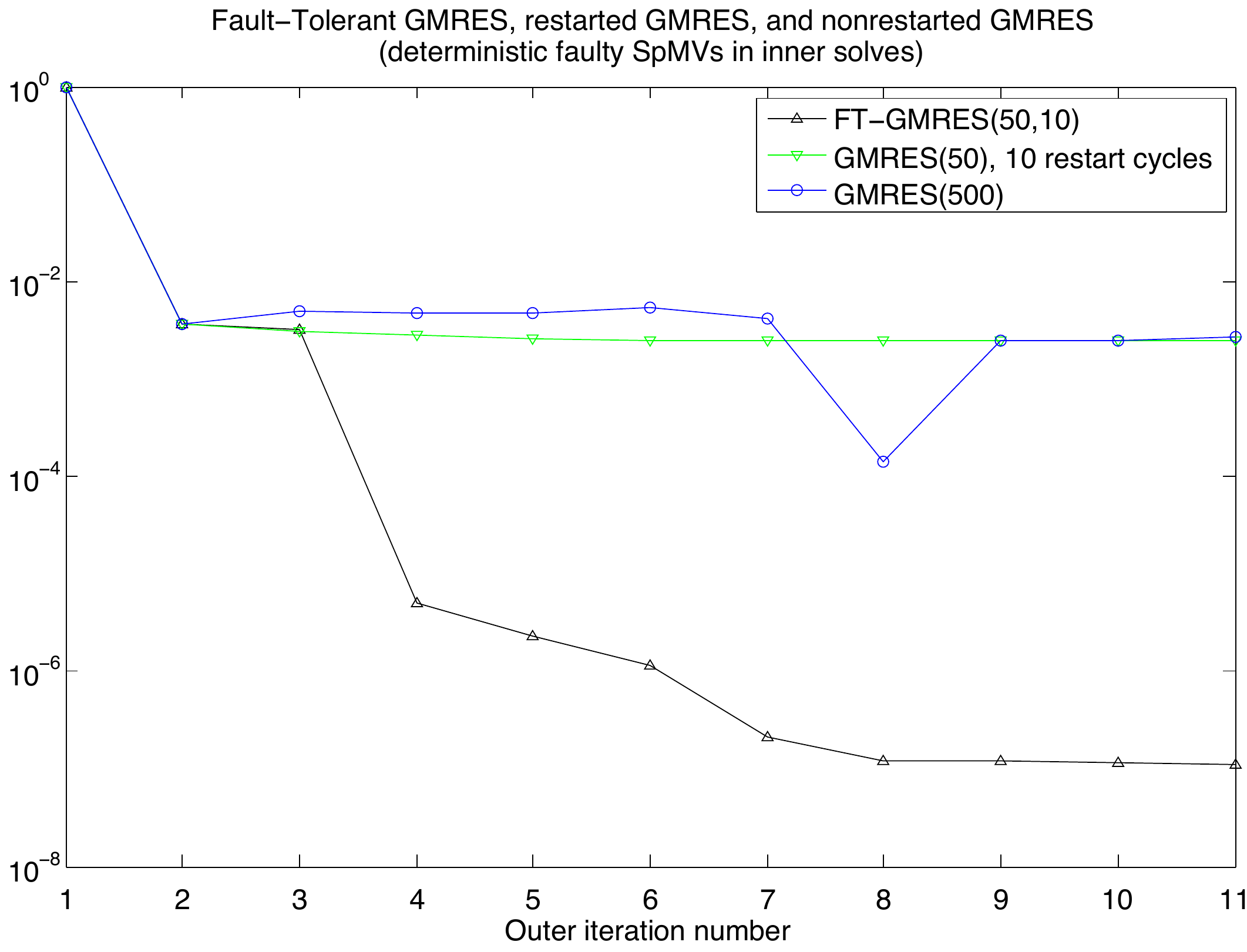}
\end{center}
\caption{\label{fig:multdcop03:restarted} FT-GMRES vs.\ GMRES on
  mult\_dcop\_03.}
\end{figure}

\begin{figure}
\begin{center}
\includegraphics[scale=0.55]{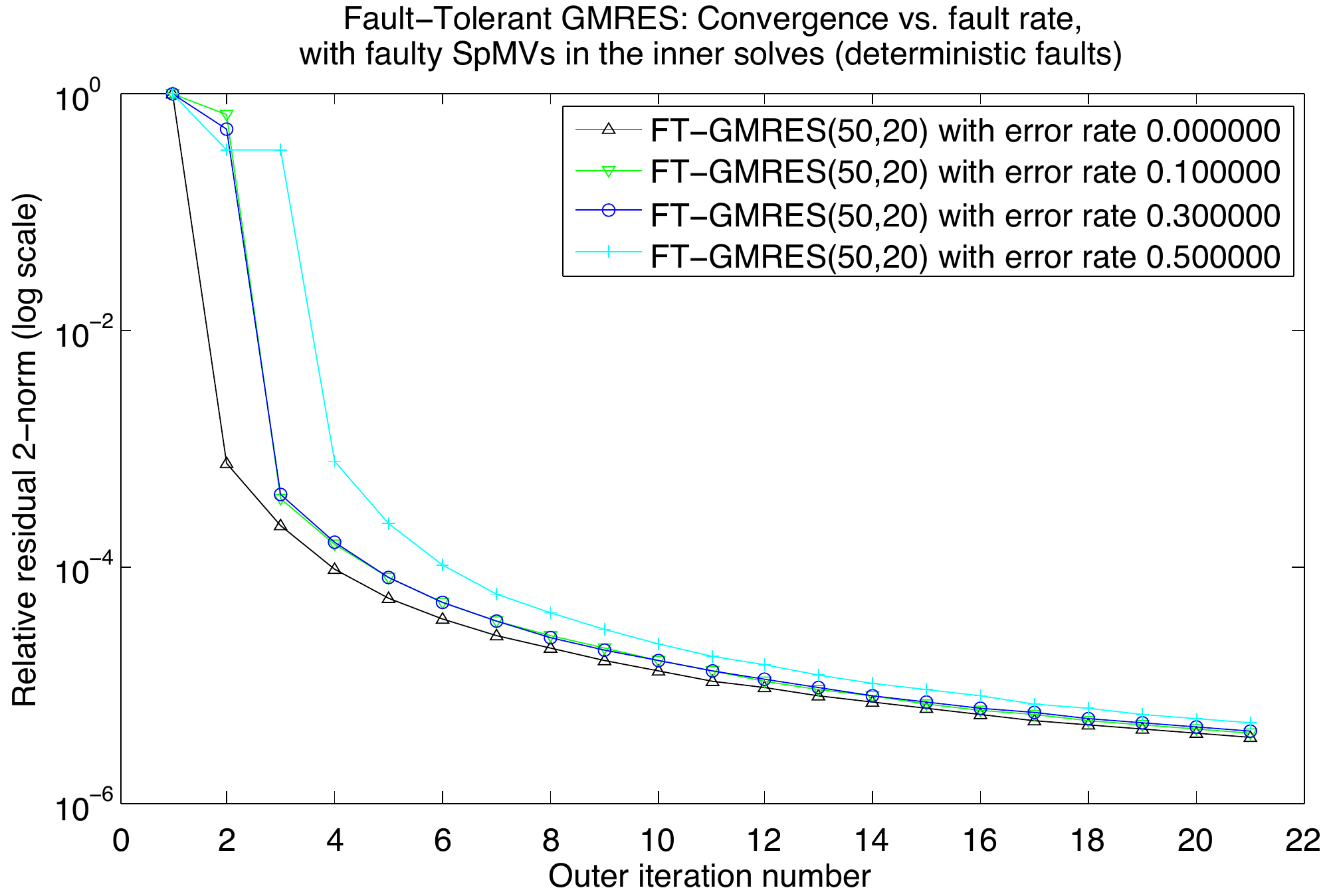}
\end{center}
\caption{\label{fig:diag10k:ftgmres} FT-GMRES on Diagonal problem,
  with different fault rates in inner solves' SpMVs.}
\end{figure}

\begin{figure}
\begin{center}
\includegraphics[scale=0.5]{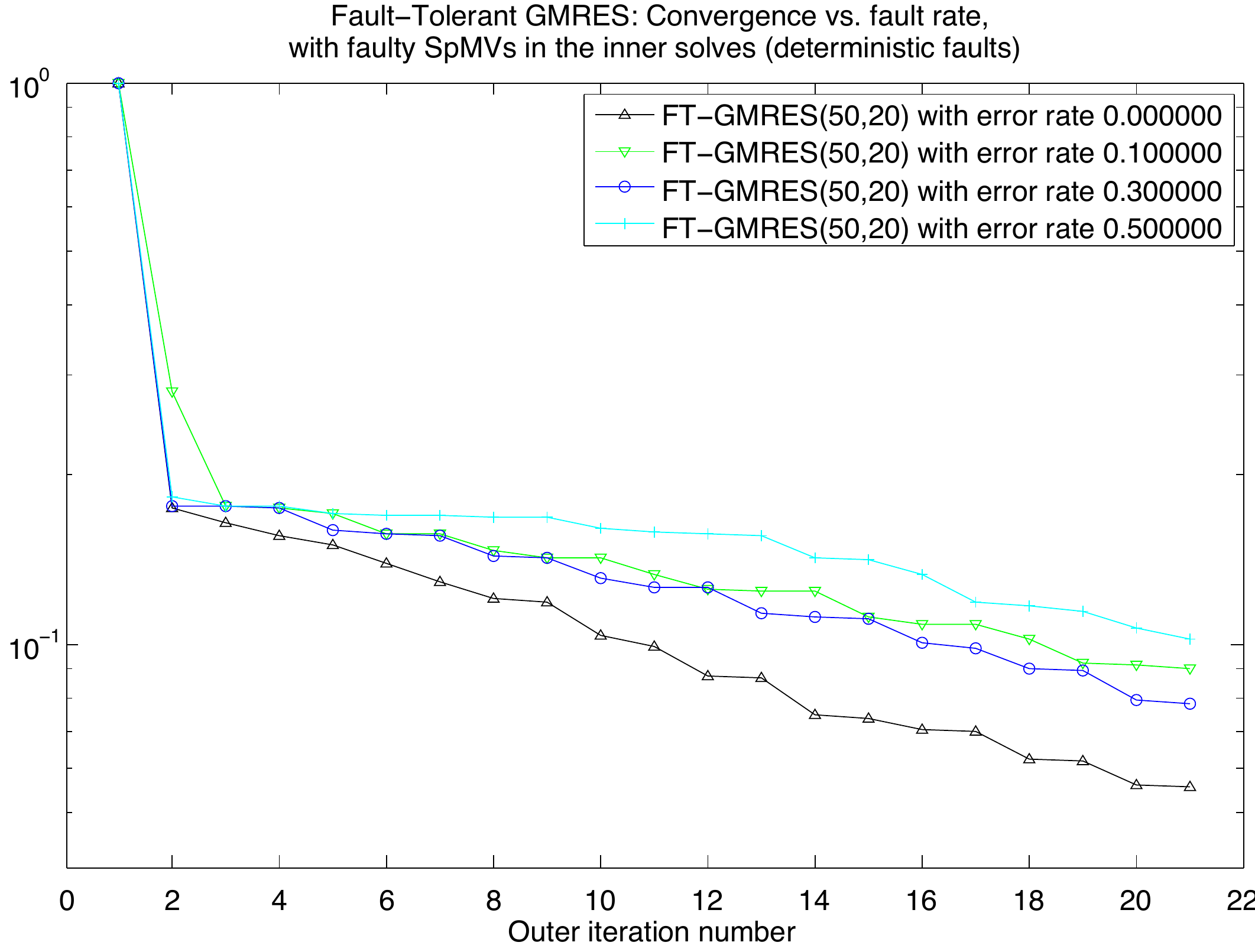}
\end{center}
\caption{\label{fig:IllStokes:ftgmres} FT-GMRES on Ill\_Stokes
  problem, with different fault rates in inner solves' SpMVs.}
\end{figure}

\begin{figure}
\begin{center}
\includegraphics[scale=0.5]{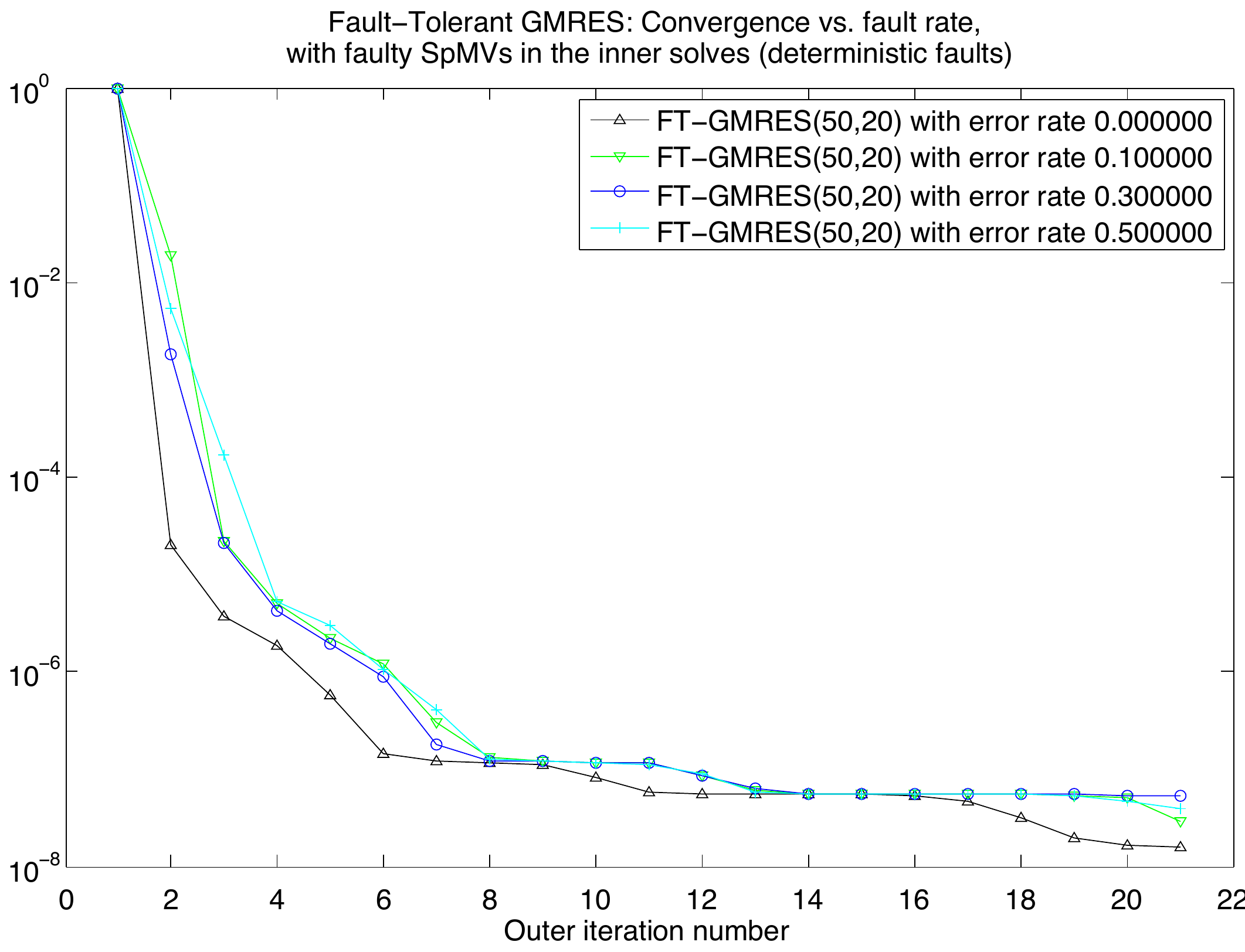}
\end{center}
\caption{\label{fig:multdcop03:ftgmres} FT-GMRES on mult\_dcop\_03
  problem, with different fault rates in inner solves' SpMVs.}
\end{figure}

\begin{figure}
\begin{center}
\includegraphics[scale=0.6]{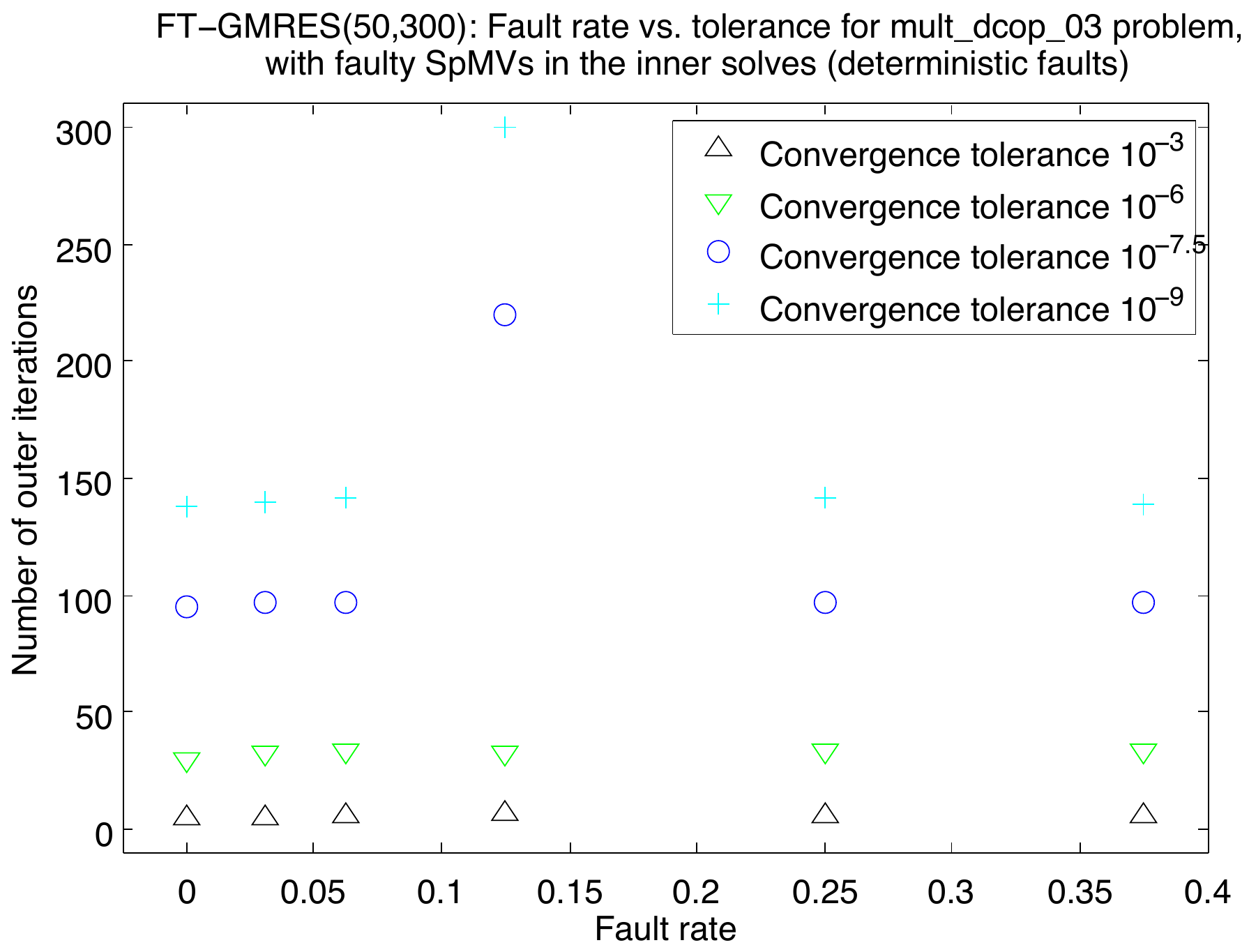}
\end{center}
\caption{\label{fig:multdcop03:table} Number of outer iterations to
  convergence for FT-GMRES (50 iterations per inner solve, max 300
  outer iterations) on mult\_dcop\_03 problem, vs.\ fault rate in the
  inner solves' SpMVs, and the outer solves' convergence tolerance.}
\end{figure}

Figures \ref{fig:diag10k:restarted}, \ref{fig:IllStokes:restarted},
and \ref{fig:multdcop03:restarted} compare FT-GMRES (50 iterations per
inner solve, 10 inner solves) with restarted GMRES (50 iterations per
restart cycle, 10 restart cycles) and nonrestarted GMRES ($500 = 50
\cdot 10$ iterations).  Every first and third out of 10 SpMVs in
GMRES, and in FT-GMRES' inner solves, are faulty.  In all cases,
FT-GMRES converges faster than the other two methods, and faults cause
restarted GMRES to stagnate or converge more slowly than FT-GMRES.
Nonrestarted GMRES' residual norm often fails to be monotonic.
Figures \ref{fig:diag10k:ftgmres}, \ref{fig:IllStokes:ftgmres}, and
\ref{fig:multdcop03:ftgmres} show only FT-GMRES (50 iterations per
inner solve, 20 inner solves), with different fault rates for SpMV
operations in the inner solves: no faults, 1 out of 10, 3 out of 10,
and 5 out of 10 SpMVs faulty.\footnote{In the 1 out of 10 case, only
  the tenth of every ten is faulty.  The 3 out of 10 case uses the
  pattern 0, 0, 0, 0, 1, 0, 0, 1, 0, 1, and the 5 out of 10 case 1, 0,
  1, 0, 1, 0, 0, 1, 0, 1.}  We found that increasing the fault rate
only decreases the FT-GMRES convergence rate gradually.  Finally,
Figure \ref{fig:multdcop03:table} shows that, barring one outlier, the
number of outer iterations to attain a given convergence rate
increases little as the fault rate increases.

%%%%%%%%%%%%%%%%%%%%%%%%%%%%%%%%%%%%%%%%%%%%%%%%%%%%%%%%%%%%%%%%%%%%%%
%%%%%%%%%%%%%%%%%%%%%%%%%%%%%%%%%%%%%%%%%%%%%%%%%%%%%%%%%%%%%%%%%%%%%%

%% file: interface.tex
%%%%%%%%%%%%%%%%%%%%%%%%%%%%%%%%%%%%%%%%%%%%%%%%%%%%%%%%%%%%%%%%%%%%%%
\section{Application / OS interface}\label{S:arch}
%%%%%%%%%%%%%%%%%%%%%%%%%%%%%%%%%%%%%%%%%%%%%%%%%%%%%%%%%%%%%%%%%%%%%%

This section describes the interface between the application and the
operating system that implements a subset of the fault detection and
recovery model described in Section \ref{S:model}.  In Section
\ref{SS:arch:design}, we present the interface itself.  Section
\ref{SS:arch:impl} outlines the implementation of this interface.
Finally, we explain in Section \ref{SS:arch:inj} our technique for
injecting artificial faults, which we use to test both the interface
and also our FT-GMRES performance prototype.  We describe fault
detection separately from fault injection, in order to emphasize that
our fault detection interface can work for actual memory faults as
well as those which the injection framework described in Section
\ref{SS:arch:inj} injects.

Our fault detection interface between the system and the application
supports both actual and artificially injected memory faults.  This
means that the FT-GMRES implementation is ready for use with existing
hardware and applications.  However, the implementation of the
system-application interface currently depends on artificial fault
injection in order to be implemented in user space on Linux.  Removing
this limitation is technically possible, but requires operating system
modifications, which prevents us from running tests on computers we do
not administer.  We leave this for future work.

%%%%%%%%%%%%%%%%%%%%%%%%%%%%%%%%%%%%%%%%%%%%%%%%%%%%%%%%%%%%%%%%%%%%%%
\subsection{Design}\label{SS:arch:design}
%%%%%%%%%%%%%%%%%%%%%%%%%%%%%%%%%%%%%%%%%%%%%%%%%%%%%%%%%%%%%%%%%%%%%%

We have designed an application / operating system (OS) interface to
support the fault and recovery models described in Section
\ref{S:model}, and implemented a library to provide this
interface. Our key design goals were to provide a simple interface for
applications and algorithmic libraries, and to support existing
OS-level interfaces to handling memory errors such as those provided
by Linux.

\begin{figure}[t!]
\begin{verbatim}
/* Register callback for handling failure in a specific 
 * allocation of failable memory at a specified byte offset 
 * and length. arg is an opaque user-supplied argument. */
typedef void (*memfail_callback_t)( void *allocation, 
                                    size_t offset, 
                                    size_t len, 
                                    void *arg);
void memfail_recover_init( memfail_callback_t cb, void *arg );

/* Mark resp. unmark memory as "failable" that was allocated
   with malloc().  Such memory should be freed with free(). */
void * malloc_failable( size_t len );
void free_failable( void *addr );
\end{verbatim}
  \caption{Application / Library interface to handle DRAM memory
    failures}
\label{fig:interface}
\end{figure}

This application level of this interface, shown in Figure
\ref{fig:interface}, focuses on marking or unmarking contiguous memory
regions that were allocated at run time using \texttt{malloc()}.  In
particular, the interface provides the application with separate calls
for marking or unmarking allocations as \emph{failable memory} -- that
is, memory in which failures will cause notifications to be sent to
the application, rather than the usual fail-stop behavior of killing
the application.  In addition, the application also registers a
callback with the library.  The callback is called once for every
active allocation when the library is notified by the OS of a detected
but uncorrected memory fault in that allocation.

In addition to this interface, we also provide a simple
producer-consumer bounded ring buffer that the application can use to
queue up a sequence of failed allocations when signaled by the
library. This ring buffer is non-blocking and atomic to allow
asynchronous callbacks from the library to enqueue failed allocations
that will be fully recovered at the end of an iteration. The
application determines the size of this buffer when it is allocated;
the number of entries needed must be sufficient to cover all of the
allocations that could plausibly fail during a single iteration. For
applications with relatively few failable allocations, this should be
a minimal number of entries.

At the OS level, the library first notifies the operating system that
it wishes to receive notifications of DRAM failures, either in general
or in specific areas of its virtual address space depending upon the
interface provided by the operating system. Second, the library keeps
track of the list of failable memory allocated by the application so
that it can call the application callback for each failed allocation
when necessary. Finally, the library handles any error notifications
from the operating system (e.g., using a Linux \texttt{SIGBUS} signal
handler) and performs OS-specific actions to clear a memory error from
a page of memory if necessary prior to notifying the application of
the error.

%%%%%%%%%%%%%%%%%%%%%%%%%%%%%%%%%%%%%%%%%%%%%%%%%%%%%%%%%%%%%%%%%%%%%%
\subsection{Implementation}\label{SS:arch:impl}
%%%%%%%%%%%%%%%%%%%%%%%%%%%%%%%%%%%%%%%%%%%%%%%%%%%%%%%%%%%%%%%%%%%%%%

We added support for handling signaled memory failures as described in
the previous section to an existing incremental checkpointing library
for Linux, libhashckpt \cite{Ferreira:11:libhashckpt}.  
% We chose this
% library because efficient incremental checkpointing requires careful
% low-level tracking of memory pages.  We were able to exploit this
% feature for a different purpose, namely the tracking of failable
% memory.  
The library also helps track application memory usage, and provides
checkpointing functionality to recover from memory failures for
applications that cannot.  Its ability to trap specific memory
accesses eases the testing of simulated memory failures, as described
later in Section~\ref{SS:arch:inj}.  Future work may include using
libhashckpt to implement efficient local save and restore of failable
data.  We modified libhashckpt to add the application API calls listed
previously in Figure~\ref{fig:interface}, with routines to mark or
unmark memory as failable.  This allocator also keeps a data structure
sorted by allocation address of failable memory allocations.

Our fault detection system assumes that when the memory controller
detects an error which ECC cannot correct, the controller notifies the
operating system using a signal that indicates which cache line in
which memory page failed.  Linux notifies the library of DRAM memory
failures, particularly failures caught by the memory scrubber using a
\texttt{SIGBUS} signal that indicates the address of the memory
\emph{page} which failed. The library then unmaps this failed page
using \texttt{munmap()}, maps in a new physical page using
\texttt{mmap()}, and calls the application-registered callback with
appropriate offset and length arguments for every failable application
allocation that overlapped with the page that included the failure.

Note that Linux currently only notifies the application of DRAM
failures detected by the memory scrubber. When the memory controller
raises an exception caused by the application attempting to consume
faulty data, Linux currently kills the faulting application. In
addition, Linux only notifies applications of the \emph{page} that
failed and expects the application to discard the entire failed
page. This approach is overly restrictive in some cases, as the
hardware notifies the kernel of the memory bus line that failed, and
some memory errors are soft and could be corrected simply by rewriting
the failed memory line.  This is not a limitation of the interface,
but of the Linux operating system itself.

%%%%%%%%%%%%%%%%%%%%%%%%%%%%%%%%%%%%%%%%%%%%%%%%%%%%%%%%%%%%%%%%%%%%%%
\subsection{Fault injection}\label{SS:arch:inj}
%%%%%%%%%%%%%%%%%%%%%%%%%%%%%%%%%%%%%%%%%%%%%%%%%%%%%%%%%%%%%%%%%%%%%%

To provide support for testing DRAM memory failures, we added support
to the incremental checkpointing library for simulating memory
failures. In particular, we added code that randomly injects errors at
a configurable rate into the application address space and uses page
protection mechanisms, i.e., \texttt{mprotect()}, to signal the
application with a \texttt{SIGSEGV} when it touches a page to which a
simulated failure has been injected. The library then catches
\texttt{SIGSEGV} and proceeds as if it had received a memory failure
on the protected page.  We model the occurrence of faults with a
Poisson distribution, with a user-specified rate (faults per MB per
hour).  We model fault locations with a uniform distribution over all
failable memory regions (i.e., those under control of our fault
detection system).  Faults that would occur in memory not currently
marked failable are not injected.

We also implemented a software simulation of a memory ``patrol
scrubber'' in the library.  The software scrubber can asynchronously
inject memory failures into the application by signaling the library
when it scrubs a memory location at which a failure has been
simulated.  For each MPI process, we start a thread for injection.
Every millisecond, that thread wakes up.  (We use the POSIX Realtime
Extension's \verb!nanosleep()! method to implement low-overhead
sleeping.)  The thread computes the number of faults that should have
been injected since the last time the thread woke up, and then injects
that many faults.  Thus, fault occurrence ``discretizes'' the Poisson
distribution, rather than obeying it exactly.  Waking up the injection
thread at longer intervals reduces its performance impact.

%% file: impl.tex
\section{Performance prototype}\label{S:impl}
%%%%%%%%%%%%%%%%%%%%%%%%%%%%%%%%%%%%%%%%%%%%%%%%%%%%%%%%%%%%%%%%%%%%%%

In this section, we describe the implementation details of our
performance prototype of FT-GMRES.  We call it a \emph{performance
  prototype} because while it may lack some features of a
production-ready implementation, we expect it to have comparable
performance.  In particular, our solver uses the same linear algebra
computational kernels and data structures as other production-ready
iterative solvers.  We made only minimal modifications to the linear
algebra objects in order to support a fault-tolerant programming
model; see the rest of this section for details.  Using the same
linear algebra objects means that when we turn off fault injection and
detection, our implementation has the same performance characteristics
as that of any other iterative solver in Trilinos.  Second, our solver
gains all the features of the linear algebra library ``for free'': in
this case, hybrid distributed- and shared-memory parallelism.  Third,
we can exploit libraries that use the linear algebra objects in order
to use our solver in a more realistic way.  For example, our
performance results in Section \ref{S:perf-exp} use a nontrivial
incomplete factorization preconditioner, and we can inject and detect
faults in the preconditioner as well.  Our solver even has the same
interface as other iterative solvers in Trilinos, so it may even be
embedded in a real application without code changes.

We call our implementation a ``prototype'' because it does not attempt
to handle all the kinds of hardware faults that we think the FT-GMRES
algorithm could handle, given the right system support.  In
particular, our system library can only intercept machine check
exceptions resulting from detected but uncorrectible DRAM faults.  Our
fault injection framework only injects faults of this kind.  In
Section \ref{SS:model:restriction}, we argue that at least for Krylov
subspace methods, bit flips are a good model for all kinds of hardware
faults, including faulty arithmetic and corrupted messages.  While our
current run-time system may not be able to handle such faults, our
more general programming model can.

%%%%%%%%%%%%%%%%%%%%%%%%%%%%%%%%%%%%%%%%%%%%%%%%%%%%%%%%%%%%%%%%%%%%%%
\subsection{Implementation outline}\label{SS:impl:org}
%%%%%%%%%%%%%%%%%%%%%%%%%%%%%%%%%%%%%%%%%%%%%%%%%%%%%%%%%%%%%%%%%%%%%%

Our prototype relies on the application / OS interface presented in
Section \ref{S:arch}.  The solver is built using components from a
slightly modified version of the Trilinos solvers framework
\cite{heroux2005overview}.  We use the implementations of GMRES
\cite{saad1986gmres} and Flexible GMRES \cite{saad1993flexible} in
Trilinos' Belos package \cite{bavier2011amesos2belos} for the inner
resp.\ outer solvers.  We preconditioned the inner iterations using
the implementation of the ILUT preconditioner \cite{saad1994ilut} in
Trilinos' Ifpack2 package \cite{williams2010ifpack2}.  All of these
algorithms use the hybrid-parallel (threads + MPI) distributed sparse
matrix and dense vector objects provided by Trilinos' Tpetra package
\cite{engineering2011baker,gradually2011baker}.  Implementing the
fault model required minor modifications to the aforementioned four
Trilinos packages, in order to demonstrate a working linear solver.
These modifications suffice to deploy FT-GMRES for use by any
application that employs Trilinos' Tpetra linear algebra stack, and
the application would only interact with them by requesting FT-GMRES
as the solver.  Our entire solver prototype required only about 3000
lines of code, not including the minor Trilinos modifications.

%%%%%%%%%%%%%%%%%%%%%%%%%%%%%%%%%%%%%%%%%%%%%%%%%%%%%%%%%%%%%%%%%%%%%%
\subsection{The Failable interface}\label{SS:impl:failable}
%%%%%%%%%%%%%%%%%%%%%%%%%%%%%%%%%%%%%%%%%%%%%%%%%%%%%%%%%%%%%%%%%%%%%%

In order for FT-GMRES to control reliability of the inner and outer
solves, we modified Trilinos to make all the large linear algebra
objects that the solver users -- the sparse matrix $A$, the
preconditioner, and vectors -- implement a \emph{Failable} abstract
interface.  The Failable interface has methods for marking, unmarking,
and checking whether the object's data are allowed to experience bit
flips.  FT-GMRES mark failability of the relevant objects on entry to
the inner solver, and unmarks them on exit.  We describe below how we
implemented this high-level interface using the low-level application
/ OS interface presented in Section \ref{S:arch}.

Trilinos is built on the Petra framework of distributed linear algebra
objects.  Petra has two implementations: Epetra (Essential Petra), and
Tpetra (Templated Petra).  We use only Tpetra for our prototype,
because Tpetra's intranode parallel support library, Kokkos
\cite{baker2010lightweight}, has the necessary features to support our
desired programming model.  In particular, Kokkos allows us to
intercept allocation and deallocation of large memory arrays, called
\emph{compute buffers}.  Linear algebra objects such as sparse
matrices, vectors, and preconditioners use compute buffers exclusively
to store data on which they plan to execute parallel kernels.  This
lets us restrict where memory faults may occur, with minimal changes
to the code of affected linear algebra objects.  Kokkos also handles
intranode parallelism in a generic way that encompasses both multicore
CPU and GPU-based hardware.  (In fact, this is why Kokkos needs
control of memory allocation; it may need to place data on a GPU or
other accelerator with a separate memory space from the CPU.)  This
lets our FT-GMRES prototype use hybrid parallelism (MPI and a
threading library of our choice) without additional effort.  Our
software prototype currently works with multiple CPU-based threading
libraries; we do not currently have GPU fault detection or injection
capability, but this could be added at the level of the application /
OS interface without changing our Trilinos modifications.

We first extended the Kokkos interface to support marking or unmarking
a compute buffer as ``failable.''  This operation directly invokes the
application / OS interface discussed in Section \ref{S:arch}.  Our
Kokkos extension gives us two ways to mark failability.  We may either
mark or unmark all subsequent allocations of compute buffers of a
particular type (e.g., \texttt{double}) as failable, or mark or unmark
a particular compute buffer.  The first option lets us experiment with
faults in Tpetra-based libraries without modifying their code.  (For
example, we can compute the sparse matrix $A$ reliably, then intercept
final assembly so that the matrix data are stored unreliably.)  The
second option -- marking each buffer individually -- lets us extend
Tpetra linear algebra objects to implement the Failable interface.

We then made Tpetra sparse matrices (CrsMatrix) and dense vectors
(MultiVector) implement the Failable interface.  Just like compute
buffers, Failable objects may be marked or unmarked failable.  Certain
data in the object may experience memory faults only if the object is
currently marked failable.  Marking a Failable object consisting of
compute buffers means marking some of its compute buffers.  The
object's implementation gets to control which compute buffers may
experience faults.  For example, our sparse matrices only mark their
entries, not the sparsity structure.  We can also compose more
complicated Failable objects out of simpler Failable objects.  For
example, an ILUT incomplete factorization preconditioner consists of
two sparse matrices (the $L$ and $U$ factors); marking the
preconditioner failable means marking the $L$ and $U$ factors
accordingly.  We modified Ifpack2's ILUT preconditioner in this way.
Finally, we added an option to Belos' GMRES solver, whether to allow
memory faults in its Krylov basis vectors.  This option was enabled
for inner solves in all of our performance experiments in Section
\ref{S:perf-exp}.  By analogy, this makes GMRES implement Failable.

%%%%%%%%%%%%%%%%%%%%%%%%%%%%%%%%%%%%%%%%%%%%%%%%%%%%%%%%%%%%%%%%%%%%%%
\subsection{Proxy for local checkpointing}\label{SS:impl:checkpointing}
%%%%%%%%%%%%%%%%%%%%%%%%%%%%%%%%%%%%%%%%%%%%%%%%%%%%%%%%%%%%%%%%%%%%%%

FT-GMRES needs to recover the correct the sparse matrix and
preconditioner after exit of each inner solve.  We implement this
using a local checkpointing scheme for failable data in the sparse
matrix and preconditioner, as we mentioned in Section \ref{S:model}.
We only save the sparse matrix values, not the sparsity structure,
since we only allow the values to experience memory faults.  Since we
do not currently have a test machine with a reliable persistent local
store, our recovery method stores the ``backup'' in nonfailable memory
as a proxy.  This serves as a performance upper bound.  We also used
this method for saving and restoring Ifpack2's ILUT preconditioner,
since it stores its factors as Tpetra sparse matrices.  This
particular preconditioner is local, so we could have just recomputed
it.  However, many preconditioners are not local, and recomputing them
often requires communication.  (The coarse-grid operators in algebraic
multigrid are a good example.)  Moreover, recomputing the ILUT
incomplete factorization each time would be expensive.

While the application / OS interface allows us to implement local
approximate repair of sparse matrix values, as described in Section
\ref{SS:model:recovery}, we have not yet implemented this technique.
The only repair we currently perform is to refresh all the sparse
matrix values from reliable backing store on exit from the inner
solve.

%%%%%%%%%%%%%%%%%%%%%%%%%%%%%%%%%%%%%%%%%%%%%%%%%%%%%%%%%%%%%%%%%%%%%%
\subsection{Operator wrappers}\label{SS:impl:wrap}
%%%%%%%%%%%%%%%%%%%%%%%%%%%%%%%%%%%%%%%%%%%%%%%%%%%%%%%%%%%%%%%%%%%%%%

In our MATLAB-based numerical experiments (see Section
\ref{SS:num-exp:framework}), we scan the vector output of each sparse
matrix or preconditioner application for invalid floating-point values
(\texttt{Inf} or \texttt{NaN}), and replace them with random data.  In
our performance prototype, we implement a similar filter.  In order to
allow thread-parallel application of the filter without a thread-safe
pseudorandom number generator, we implement ``repairing'' the vector
by replacing an invalid value with the arithmetic mean of its nearest
valid neighbors (within a fixed-size window of neighbors).  We use the
Kokkos parallel framework to apply this kernel in parallel.

%%% Local Variables: 
%%% mode: latex
%%% TeX-master: "paper"
%%% End: 

%% file: experiments.tex
%%%%%%%%%%%%%%%%%%%%%%%%%%%%%%%%%%%%%%%%%%%%%%%%%%%%%%%%%%%%%%%%%%%%%%
\section{Performance experiments}\label{S:perf-exp}
%%%%%%%%%%%%%%%%%%%%%%%%%%%%%%%%%%%%%%%%%%%%%%%%%%%%%%%%%%%%%%%%%%%%%%

In this section, we present the results of experiments with our
FT-GMRES performance prototype.  We demonstrate the solver working
with multiple threads and MPI processes, using production-quality
software components from the Trilinos framework.  
% Finally, we show the
% performance overhead of fault detection and injection.

We tested FT-GMRES using the development (10.7) version of Trilinos,
on a Intel Xeon X5570 (8 cores, 2.93 GHz) CPU with 12 GB of main
memory.  We chose for our test matrix an ill-conditioned Stokes
partial differential equation discretization
\texttt{Szczerba/Ill\_Stokes} from the University of Florida Sparse
Matrix Collection (UFSMC) \cite{UFSMC}.  It has 25,187 rows and
columns, 193,216 stored entries, and an estimated 1-norm condition
number of $4.85 \times 10^9$.  For initial experiments, we chose a
uniform $[-1,1]$ pseudorandom right-hand side.

\begin{figure}
\begin{center}
\includegraphics[scale=0.8]{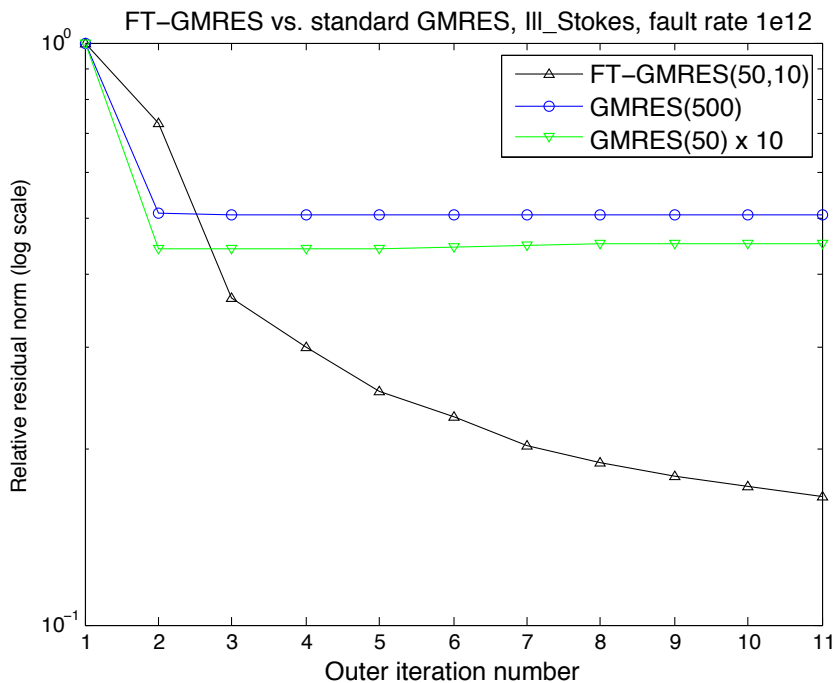}
\end{center}
\caption{\label{fig:faultRate1e12} FT-GMRES (10 outer iterations, 50
  inner iterations each), 500 iterations of non-restarted GMRES, and 10
  restart cycles (50 iterations each) of restarted GMRES.  (Down is
  good.)}
\end{figure}

We ran FT-GMRES with 10 outer iterations.  Each inner solve used 50
iterations of standard GMRES (without restarting),
right-preconditioned by ILUT (see e.g., Saad \cite{saad2003iterative})
with level 2 fill, zero drop tolerance, 1.5 relative threshold and 0.1
absolute threshold.  (These are not necessarily reasonable ILUT
parameters, but they ensure a valid preconditioner for the problem
tested.)  We compared FT-GMRES with standard GMRES both with and
without restarting: 500 iterations of each, restarting if applicable
every 50 iterations.  (This makes the memory usage of the two methods
approximately comparable.)  We set no convergence criteria except for
iteration counts, so that we could fully observe the behavior of the
methods.  Our initial experiments use random fault injection at a rate
of 1000 faults per megabyte per hour, which is high but demonstrates
the solver's fault-tolerance capabilities.  Faults were allowed to
occur in floating-point data belonging to the matrix and the ILUT
preconditioner.  Furthermore, to demonstrate the value of algorithmic
approaches, our restarted GMRES implementation imitated FT-GMRES by
also refreshing the matrix and ILUT preconditioner from reliable
storage before every restart cycle.  (We optimized by not refreshing
if no memory faults were detected.)

Figure \ref{fig:faultRate1e12} shows our convergence results.
FT-GMRES' reliable outer iteration makes it able to roll forward
through faults and continue convergence.  The fault-detection
capabilities discussed earlier in this work let FT-GMRES refresh
unreliable data only when necessary, so that memory faults appear
transient to the solver.

%%% Local Variables: 
%%% mode: latex
%%% TeX-master: "paper"
%%% End: 

%% file: conclusion.tex
%%%%%%%%%%%%%%%%%%%%%%%%%%%%%%%%%%%%%%%%%%%%%%%%%%%%%%%%%%%%%%%%%%%%%%
\section{Conclusion}\label{S:conclusion}
%%%%%%%%%%%%%%%%%%%%%%%%%%%%%%%%%%%%%%%%%%%%%%%%%%%%%%%%%%%%%%%%%%%%%%

In this paper, we showed that our fault-tolerant iterative linear
solver, FT-GMRES, can converge despite memory faults, whereas other
iterative solvers could not.  We demonstrated this both with numerical
experiments using a Matlab prototype, and with a performance prototype
that uses realistic fault injection and detection, implemented using
production-grade solver components.

Our experiments also show that FT-GMRES' convergence rate degrades
gradually as the fault rate is increased, and that increasing the
fault rate only modestly increases the total number of iterations (and
therefore the total cost).  While more experiments are needed, we
think FT-GMRES and fault-tolerant iterative methods in general have
great potential to improve solver robustness and relax hardware
reliability constraints.  The basic approaches we have used can be
applied to many algorithms, greatly reducing the impact of the soft
faults that are expected on future computing systems.

Our work has also opened up interesting collaborations with systems
researchers, to develop programming interfaces for varying
reliability, reporting faults, and selective checkpointing.  These
collaborations have the potential to influence hardware-software
codesign, especially at extreme scales, where energy requirements will
force system designers to reduce hardware reliability and rely more on
software approaches.  Fault-tolerant algorithms thus have the
potential to influence computer hardware in a way analogous to RISC
(Reduced Instruction Set Computer) architectures
\cite{patterson1981risc}, by encouraging beneficial trade-offs between
hardware and software.

%%% Local Variables: 
%%% mode: latex
%%% TeX-master: "paper"
%%% End: 

%% file: ack.tex
\section*{Acknowledgment}\label{S:ack}
%%%%%%%%%%%%%%%%%%%%%%%%%%%%%%%%%%%%%%%%%%%%%%%%%%%%%%%%%%%%%%%%%%%%%%

%Patrick G.~Bridges: 
This work was supported in part by a faculty sabbatical appointment
from Sandia National Laboratories and a grant from the U.S. Department
of Energy and a grant from DOE Office of Science, Advanced Scientific
Computing research, under award number DE-SC0005050, program manager
Sonia Sachs.

Sandia National Laboratories is a multiprogram laboratory managed and
operated by Sandia Corporation, a wholly owned subsidiary of Lockheed
Martin Corporation, for the U.S.\ Department of Energy's National
Nuclear Security Administration under contract DE-AC04-94AL85000.

%% file: paper.bbl
\begin{thebibliography}{10}

\bibitem{asanovic2006landscape}
{\sc K.~Asanovic, R.~Bodik, B.~C. Catanzaro, J.~J. Gebis, P.~Husbands,
  K.~Keutzer, D.~A. Patterson, W.~L. Plishker, J.~Shalf, S.~W. Williams, and
  K.~A. Yelick}, {\em The {L}andscape of {P}arallel {C}omputing {R}esearch: {A}
  {V}iew from {B}erkeley}, Tech. Rep. UCB/EECS-2006-183, EECS Department,
  University of California, Berkeley, Dec 2006.
\newblock See also the journal article Asanovic et al.\
  \cite{asanovic2009landscape}.

\bibitem{asanovic2009landscape}
{\sc K.~Asanovic, R.~Bodik, J.~W. Demmel, T.~Keaveny, K.~Keutzer,
  J.~Kubiatowicz, N.~Morgan, D.~A. Patterson, K.~Sen, J.~Wawrzynek, D.~Wessel,
  and K.~A. Yelick}, {\em A {V}iew of the {P}arallel {C}omputing {L}andscape},
  Communications of the ACM, 52 (2009), pp.~56--67.

\bibitem{bahi2007parallel}
{\sc J.~Bahi, S.~Contasset-Vivier, and R.~Couturier}, {\em Parallel Iterative
  Algorithms: From Sequential to Grid Computing}, Chapman and Hall / CRC, 2007.

\bibitem{engineering2011baker}
{\sc C.~G. Baker}, {\em Engineering a kernel-agnostic distributed linear
  algebra library for multi/many-core in {T}rilinos}.
\newblock Refereed abstract in minisymposium ``Parallel Programming Models and
  Algorithms for Scalable Manycore Systems,'' SIAM Conference on Computational
  Science and Engineering 2011, Reno, NV, March 2011.

\bibitem{gradually2011baker}
\leavevmode\vrule height 2pt depth -1.6pt width 23pt, {\em Gradually
  transitioning library users to the hybrid-parallel paradigm}.
\newblock Refereed abstract in minisymposium ``Creating the Next Generation of
  High Performance Numerical Computing Capabilities,'' International Congress
  on Industrial and Applied Mathematics (ICIAM) 2011, Vancouver, B.C., July
  2011.

\bibitem{baker2010lightweight}
{\sc C.~G. Baker, H.~C. Edwards, M.~A. Heroux, and A.~B. Williams}, {\em A
  light-weight {API} for portable multicore programming}, in The 18th Euromicro
  International Conference on Parallel, Distributed and Network-Based Computing
  (PDP 2010), February 2010.

\bibitem{bavier2011amesos2belos}
{\sc E.~Bavier, M.~Hoemmen, S.~Rajamanickam, and H.~Thornquist}, {\em {Amesos2}
  and {Belos}: Direct and iterative solvers for large sparse linear systems},
  2011.
\newblock Submitted to Scientific Programming.

\bibitem{blackford1996practical}
{\sc L.~S. Blackford, A.~Cleary, J.~Demmel, I.~Dhillon, J.~Dongarra,
  S.~Hammarling, A.~Petitet, H.~Ren, K.~Stanley, and R.~C. Whaley}, {\em
  Practical experience in the dangers of heterogeneous computing}, Tech. Rep.
  UT-CS-96-330, University of Tennessee, Knoxville, July 1996.
\newblock LAPACK Working Note \#112.

\bibitem{blackford1997practical}
\leavevmode\vrule height 2pt depth -1.6pt width 23pt, {\em Practical experience
  in the numerical dangers of heterogeneous computing}, ACM Trans. Math.
  Softw., 23 (1997), pp.~133--147.

\bibitem{Bronevetsky:2008:Soft}
{\sc G.~Bronevetsky and B.~de~Supinski}, {\em Soft error vulnerability of
  iterative linear algebra methods}, in Proceedings of the 22nd Annual
  International Conference on Supercomputing, ICS '08, New York, NY, USA, 2008,
  ACM, pp.~155--164.

\bibitem{buttari2006computations}
{\sc A.~Buttari, J.~Dongarra, J.~Kurzak, P.~Luszczek, and S.~Tomov}, {\em
  Computations to enhance the performance while achieving the 64-bit accuracy},
  Tech. Rep. UT-CS-06-584, University of Tennessee Knoxville, November 2006.
\newblock LAPACK Working Note \#180.

\bibitem{castro2002practical}
{\sc M.~Castro and B.~Liskov}, {\em Practical {B}yzantine fault tolerance and
  proactive recovery}, ACM Transactions on Computer Systems, 20 (2002),
  pp.~398--461.

\bibitem{chen2005compiler}
{\sc G.~Chen, M.~Kandemir, M.~J. Irwin, and G.~Memik}, {\em Compiler-directed
  selective data protection against soft errors}, in In Proceedings of the Asia
  and South Pacific Design Automation Conference (ASPDAC), Shanghai, China,
  January 2005.

\bibitem{Chishti:2009:ICL:1669112.1669126}
{\sc Z.~Chishti, A.~R. Alameldeen, C.~Wilkerson, W.~Wu, and S.-L. Lu}, {\em
  Improving cache lifetime reliability at ultra-low voltages}, in Proceedings
  of the 42nd Annual IEEE/ACM International Symposium on Microarchitecture,
  MICRO 42, New York, NY, USA, 2009, ACM, pp.~89--99.

\bibitem{UFSMC}
{\sc T.~A. Davis and Y.~Hu}, {\em The {U}niversity of {F}lorida {S}parse
  {M}atrix {C}ollection}, ACM Trans. Math. Softw. (to appear),  (2011).

\bibitem{Ferreira:11:libhashckpt}
{\sc K.~B. Ferreira, R.~Riesen, R.~Brightwell, P.~G. Bridges, and D.~Arnold},
  {\em Libhashckpt: Hash-based incremental checkpointing using {GPUs}}, in
  Proceedings of the 18th EuroMPI Conference, Santorini, Greece, September
  2011.

\bibitem{golub1999inexact}
{\sc G.~H. Golub and Q.~Ye}, {\em Inexact preconditioned conjugate gradient
  method with inner-outer iteration}, SIAM J. Sci. Comput., 21 (1999),
  pp.~1305--1320.

\bibitem{greenbaum1996any}
{\sc A.~Greenbaum, V.~Ptak, and Z.~Strakos}, {\em Any nonincreasing convergence
  curve is possible for {GMRES}}, SIAM J. Matrix Anal. Appl., 17 (1996),
  pp.~465--469.

\bibitem{Haque:2010:HDS:1844765.1845231}
{\sc I.~S. Haque and V.~S. Pande}, {\em Hard data on soft errors: A large-scale
  assessment of real-world error rates in {GPGPU}}, in Proceedings of the 2010
  10th IEEE/ACM International Conference on Cluster, Cloud and Grid Computing,
  CCGRID '10, Washington, DC, USA, 2010, IEEE Computer Society, pp.~691--696.

\bibitem{heroux2005overview}
{\sc M.~A. Heroux, R.~A. Bartlett, V.~E. Howle, R.~J. Hoekstra, J.~J. Hu, T.~G.
  Kolda, R.~B. Lehoucq, K.~R. Long, R.~P. Pawlowski, E.~T. Phipps, A.~G.
  Salinger, H.~K. Thornquist, R.~S. Tuminaro, J.~M. Willenbring, A.~Williams,
  and K.~S. Stanley}, {\em An overview of the {T}rilinos project}, ACM Trans.
  Math. Softw., 31 (2005), pp.~397--423.

\bibitem{hoemmen2011fault}
{\sc M.~A. Heroux and M.~Hoemmen}, {\em Fault-tolerant iterative methods via
  selective reliability}, Tech. Rep. SAND2011-3915 C, Sandia National
  Laboratories, 2011.
\newblock Available at \url{http://www.sandia.gov/~maherou/}.

\bibitem{howle2010soft:Copper}
{\sc V.~E. Howle}, {\em Soft errors in linear solvers as integrated components
  of a simulation}, in Presented at the Copper Mountain Conference on Iterative
  Methods, Copper Mountain, CO, April 9, 2010.

\bibitem{huang1984algorithm}
{\sc K.-H. Huang and J.~A. Abraham}, {\em Algorithm-based fault tolerance for
  matrix operations}, IEEE Transactions on Computers, C-33 (1984).

\bibitem{Karnik:2004:CSE:1032295.1032595}
{\sc T.~Karnik, P.~Hazucha, and J.~Patel}, {\em Characterization of soft errors
  caused by single event upsets in {CMOS} processes}, IEEE Trans. Dependable
  Secur. Comput., 1 (2004), pp.~128--143.

\bibitem{kogge:exa}
{\sc P.~M.~e. Kogge}, {\em {ExaScale Computing Study: Technology Challenges in
  Achieving Exascale Systems}}, tech. rep., University of Notre Dame CSE
  Department Technical Report, TR-2008-13, September 28, 2008.

\bibitem{lamport1982byzantine}
{\sc L.~Lamport, R.~Shostak, and M.~Pease}, {\em The {B}yzantine {G}enerals
  {P}roblem}, ACM Transactions on Programming Languages and Systems, 4 (1982),
  pp.~382--401.

\bibitem{Langou:2007:RPI:1350656.1350657}
{\sc J.~Langou, Z.~Chen, G.~Bosilca, and J.~Dongarra}, {\em Recovery patterns
  for iterative methods in a parallel unstable environment}, SIAM J. Sci.
  Comput., 30 (2007), pp.~102--116.

\bibitem{Lin:2010:TLM:1838773.1839098}
{\sc P.~T. Lin and J.~N. Shadid}, {\em Towards large-scale multi-socket,
  multicore parallel simulations: Performance of an {MPI}-only semiconductor
  device simulator}, J. Comput. Phys., 229 (2010), pp.~6804--6818.

\bibitem{malkowski2010analyzing}
{\sc K.~Malkowski, P.~Raghavan, and M.~Kandemir}, {\em Analyzing the soft error
  resilience of linear solvers on multicore processors}, in Proceedings of the
  Twenty-Fourth IEEE International Parallel and Distributed Processing
  Symposium (IPDPS), 2010.

\bibitem{Miskov-Zivanov:2007:SER:1266366.1266680}
{\sc N.~Miskov-Zivanov and D.~Marculescu}, {\em Soft error rate analysis for
  sequential circuits}, in Proceedings of the Conference on Design, Automation
  and Test in Europe, DATE '07, San Jose, CA, USA, 2007, EDA Consortium,
  pp.~1436--1441.

\bibitem{parhami1997defect}
{\sc B.~Parhami}, {\em Defect, fault, error, {\dots}, or failure?}, IEEE
  Transactions on Reliability, 46 (1997).

\bibitem{patterson1981risc}
{\sc D.~A. Patterson and C.~H. Sequin}, {\em {RISC} {I}: A {R}educed
  {I}nstruction {S}et {VLSI} {C}omputer}, in Proceedings of the Eighth Annual
  Symposium on Computer Architecture (ISCA), IEEE Computer Society Press, 1981.

\bibitem{rosenblum2004reincarnation}
{\sc M.~Rosenblum}, {\em The reincarnation of virtual machines}, Queue, 2
  (2004).

\bibitem{saad1993flexible}
{\sc Y.~Saad}, {\em A flexible inner-outer preconditioned {GMRES} algorithm},
  SIAM J. Sci. Comput., 14 (1993), pp.~461--469.

\bibitem{saad1994ilut}
\leavevmode\vrule height 2pt depth -1.6pt width 23pt, {\em {ILUT}: a dual
  threshold incomplete {ILU} factorization}, Num. Lin. Alg. Appl., 1 (1994),
  pp.~387--402.

\bibitem{saad2003iterative}
\leavevmode\vrule height 2pt depth -1.6pt width 23pt, {\em Iterative Methods
  for Sparse Linear Systems}, SIAM, Philadelphia, second~ed., 2003.

\bibitem{saad1986gmres}
{\sc Y.~Saad and M.~H. Schultz}, {\em {GMRES}: A generalized minimal residual
  algorithm for solving nonsymmetric linear systems}, SIAM J. Sci. Statist.
  Comput., 7 (1986), pp.~856--869.

\bibitem{schroeder2009dram}
{\sc B.~Schroeder, E.~Pinheiro, and W.-D. Weber}, {\em {DRAM} errors in the
  wild: A large-scale field study}, in SIGMETRICS / Performance 2009, June
  15--19, 2009, Seattle, WA, USA, 2009.

\bibitem{simonici2003flexible}
{\sc V.~Simonici and D.~B. Szyld}, {\em Flexible inner-outer {K}rylov subspace
  methods}, SIAM J. Numer. Anal., 40 (2003), pp.~2219--2239.

\bibitem{simonici2003theory}
\leavevmode\vrule height 2pt depth -1.6pt width 23pt, {\em Theory of inexact
  {K}rylov subspace methods and applications to scientific computing}, SIAM J.
  Sci. Comput., 25 (2003), pp.~454--477.

\bibitem{smith2005architecture}
{\sc J.~E. Smith and R.~Nair}, {\em The architecture of virtual machines},
  Computer, 38 (2005), pp.~32--38.

\bibitem{stewart1993updating}
{\sc G.~W. Stewart}, {\em Updating a rank-revealing {$ULV$} decomposition},
  SIAM J. Matrix Anal. Appl., 14 (1993), pp.~494--499.

\bibitem{szyld2001fqmr}
{\sc D.~B. Szyld and J.~A. Vogel}, {\em {FQMR}: A flexible quasi-minimal
  residual method with inexact preconditioning}, SIAM J. Sci. Comput., 23
  (2001), pp.~363--380.

\bibitem{eshof2004inexact}
{\sc J.~van~den Eshof and G.~L.~G. Sleijpen}, {\em Inexact {K}rylov subspace
  methods for linear systems}, SIAM J. Matrix Anal. Appl., 26 (2004),
  pp.~125--153.

\bibitem{wilkinson1963rounding}
{\sc J.~H. Wilkinson}, {\em Rounding Errors in Algebraic Processes}, Prentice
  Hall, Englewood Cliffs, NJ, 1963.
\newblock This work has been republished in an inexpensive Dover paperback
  edition.

\bibitem{williams2010ifpack2}
{\sc A.~Williams}, {\em New {T}rilinos package: {Ifpack2}}, 2010.
\newblock Presentation at Trilinos Users' Group meeting, Sandia National
  Laboratories, Albuquerque, New Mexico, 04 November 2010. Video available at
  \url{http://trilinos.sandia.gov/events/trilinos_user_group_2010/}.

\bibitem{Yang:2005:LRS:1120725.1120957}
{\sc S.~Yang, W.~Wolf, W.~Wang, N.~Vijaykrishnan, and Y.~Xie}, {\em Low-leakage
  robust {SRAM} cell design for sub-100nm technologies}, in Proceedings of the
  2005 Asia and South Pacific Design Automation Conference, ASP-DAC '05, New
  York, NY, USA, 2005, ACM, pp.~539--544.

\end{thebibliography}
